\documentclass[a4,11pt]{article}
\usepackage[dvipdfmx]{graphicx}
\usepackage{lscape}
\usepackage{enumerate}
\usepackage{latexsym}
\usepackage{amsthm}
\usepackage{subfig}
\usepackage{color}
\usepackage{array}

\usepackage{framed}

\newtheorem{theorem}{Theorem}

\newtheorem{lemma}{Lemma}

\usepackage{amsmath,amssymb}

\DeclareMathOperator{\argmin}{argmin}

\def\vec#1{\mbox{\boldmath $#1$}}
\begin{document}

\title{Single-index models for extreme value index regression}

\author{
{\sc Takuma Yoshida}$^{1}$\\
$^{1}${\it Kagoshima University, Kagoshima 890-8580, Japan}\\
{\it E-mail: yoshida@sci.kagoshima-u.ac.jp}
}

\date{\empty}
\maketitle

\begin{abstract}
Since the extreme value index (EVI) controls the tail behavior of the
distribution function, the estimation of EVI is a very important topic in
extreme value theory. Recent contributions have focused on nonparametric regression approaches with covariates for the estimation of
EVI. However, for high-dimensional settings, the fully nonparametric
estimator faces the curse of dimensionality. To resolve this, we apply the
single index model to EVI regression under a Pareto-type tailed
distribution. We study the penalized maximum likelihood estimation of
the single index model. The asymptotic properties of the estimator are
also developed. Numerical studies are presented to demonstrate the
efficiency of the proposed model.
\end{abstract}

{\it Keywords:
Extreme value index; Heavy tail; Pareto-type model; Peak over threshold; Penalized spline; Single index model}

{\it MSC codes: 62G08, 62G20, 62G32}

\section{Introduction}

Analyzing the probability of occurrence of a rare event is important to evaluate risk assessment in diverse fields such as meteorology, economics, sociology, ecology, and life sciences. 
A rare event is defined as one for which data values are extremely high or low or which are located at the tail of the distribution. 
Extreme value theory (EVT) is an efficient statistical tool for investigating the tail behavior of a distribution. 
Many authors have developed the method, theory, and application of EVT, as summarized by Beirlant et al. (2004), de Haan and Ferreira (2006), and Dey and Yan (2016). 
The tail behavior of the distribution is classifiable into three types: heavy tail, light tail, and short tail. 
Such division is controlled by a parameter called the extreme value index (EVI). 
The distribution has a heavy tail if EVI is positive. 
The light and short tails respectively correspond to zero and negative EVI. 
As described in this paper, we specifically examine the case of heavy tail or positive EVI because estimation of the positive EVI is more difficult than it is in other cases. 
The Hill estimator, proposed by Hill (1975), is known as the fundamental estimator for positive EVI.
As described in this paper, we specifically examine the case of positive EVI and assume a Pareto-type distribution.
The Hill estimator is related closely to the maximum likelihood estimator for the EVI under the Pareto-type distribution.
As explained below, our proposed estimator can be regarded as a covariate-dependent extension of the Hill estimator.

In recent years, rapid development has occurred in the estimation of the conditional EVI with covariate information in the context of regression. 
The nonparametric estimator of conditional EVI was suggested by Gardes (2010), Stupfler (2013), Daouia et al. (2013), Stupfler and Gardes (2014), Goegebeur et al. (2014), Goebegeur et al. (2015), and Ma et al. (2020).
However, in a high-dimensional setting, such an estimator presents the difficulty of the curse of dimensionality. For that reason, the efficiency of a fully nonparametric estimator cannot be guaranteed. 
Therefore, for high-dimensional covariates, one must adopt flexible modelling of the target function to avoid the curse of dimensionality. 
The linear model, proposed by Wang and Tsai (2009), is the classical approach used for flexible modelling of the target function. 
However, the linear model is unable to capture the behavior of the data having a nonlinear structure. 
In other words, the linear model is too restrictive to account for the complexity of the data.
As a flexible semiparametric approach, Chavez-Demoulin and Davison (2005) and Youngman (2019) have used the generalized additive model.
Li et al. (2020) conducted EVI regression with a partially linear model. 
Wang and Li (2012) and Wang et al. (2013) studied the conditional Hill estimator from linear extremal quantile regression.
The varying coefficient model was developed by Ma et al. (2019) and Momoki and Yoshida (2024). 
Although the study of the EVI modelling has persisted as a topic of interest in recent years, it has never been considered in relation to the single index model, although the single index model is also flexible modelling approach.
This gap in the related literature motivates us to introduce the single index model in EVI estimation with a large dimension of covariates.

Ichimura (1993) and H$\ddot{a}$rdle et al. (1993) proposed the single index model in mean regression. 
This model, known as the semiparametric model, is structured as a hybrid of the linear transformation of covariates and one-dimensional nonparametric function. 
Hall (1989), Horowits and H$\ddot{a}$rdle (1993), Carrol et al. (1997), Yu and Ruppert (2002), Wang and Yang (2009), and Kuchibhotla and Patra (2016) have developed the single index model in mean regression. 
In quantile regression, Wu et al. (2010), Zu et al. (2012), and Ma and He (2015) have studied the single index model. 
Gardes (2018) and Xu et al. (2022) considered the usage of the single index model in extremal quantile regression. 
That earlier work has motivated us to apply the single index model to EVI regression with large dimensional covariates. 
From the method presented by Gardes (2018) and by Xu et al. (2022), EVI can be estimated as the conditional Hill estimator using a conditional quantile with several quantile levels. 
However, the method presented by Gardes (2018) is complicated. Moreover, it entails high computational cost. 
Furthermore, the single index parameter depends on the quantile level. Therefore, the obtained EVI estimator has no the single index structure. 
Xu et al. (2022) assumed the linear model as the conditional quantile. 
However, for the tail quantile, the linearity assumption is too restrictive for more general settings. 
Bousebata et al. (2023) and Aghbalou et al. (2024) also consider the single-index or multi-index structure in extreme value analysis, but they do not directly and specifically examine estimation of the EVI function. 
Unlike earlier studies, our goal is estimation of the single-index parameter and EVI function simultaneously. 

In single index models, it is necessary to estimate the linear coefficient parameter vector and the one-dimensional nonlinear function. 
First, we assume that the Pareto-type-tailed model is the conditional distribution of the response variable as a function of the covariates. 
Then, the single index parameter and nonlinear function including EVI is estimated via the maximum likelihood method after choosing extreme data using the peak over threshold (POT) method. 
We estimate the nonlinear component of the single-index model using penalized splines, a standard approach in semi-parametric modeling.
We study the asymptotic distribution and the rate of convergence of the proposed estimator. 
Based on these results, we can verify whether the proposed single index model overcomes the difficulty of the curse of dimensionality. 
The finite sample performance of the proposed single index model is examined using a Monte Carlo simulation. 
We also report an empirical data example using motor bike insurance data presented by Ohlsson and Johansson (2010). 

Next, we explain why the spline method is used instead of other methods, such as the kernel smoothers, for estimating the nonlinear part.
According to Yu and Ruppert (2002) and Wang and Yang (2009), the spline method is computationally more efficient than the kernel smoothers in the single index model. 
Furthermore, from a recent study of the regression with extreme value analysis, Youngman (2022) has proposed a very useful R-package called {\sf evgam}. 
The smoothing method used in {\sf evgam} is mainly splines. 
Consequently, the demand for the methodology and the theory of the spline method is expected to increase in the field of EVT.  
This expected demand motivates us to examine the spline method specifically in this study.

The remainder of the paper is organized as presented below.
Section 2 sets the single index model for EVI regression, the estimation procedure of the maximum likelihood method and tuning parameter selection. 
Asymptotic theory for the proposed estimator is established later in Section 3. 
The simulation study is described in Section 4. The empirical data example is given in Section 5.
Thereafter, Section 6 concludes the paper. 
As presented in the Appendix, we review the important properties of splines and the technical lemmas and the proof of theorems are also provided. 

\section{Single Index Model}

\subsection{Model setting}

Consider the random pair $(Y,\vec{X})$ with the response $Y\in\mathbb{R}_+$ and the covariate $\vec{X}=(X_1,\ldots,X_p)\in{\cal X}\subset\mathbb{R}^p$. 
As described in this paper, the domain of covariate ${\cal X}$ is compact space.
Let $F(y|\vec{x})=P(Y\leq y|\vec{X}=\vec{x})$ be the conditional distribution function of $Y$ given $\vec{X}=\vec{x}=(x_1,\ldots,x_p)$. 
We then assume that $Y$ given $\vec{X}=\vec{x}$ is distributed as the class of Pareto-type tailed distributions defined as
\begin{eqnarray}
P(Y>y|\vec{X}=\vec{x})=1-F(y|\vec{x})= y^{ -1/\gamma^*(\vec{x})}L(y|\vec{x}), \label{paretotail}
\end{eqnarray}
where $\gamma^*(\vec{x})>0$ is the EVI function and $L$ is a slowly varying function satisfying 
\[
\lim_{y\rightarrow \infty} L(ay|\vec{x})/L(y|\vec{x})=1
\] for all $\vec{x}\in{\cal X}$ and $a>0.$
As described in this paper, the slowly varying function $L$ is assumed to belong to the Hall class (Hall 1982) as
\begin{eqnarray}
L(y|\vec{x})= \ell_0(\vec{x})+ \ell_1(\vec{x})y^{-\beta(\vec{x})}+\nu(y|\vec{x}), \label{slowly}
\end{eqnarray}
where for any $\vec{x}\in{\cal X}$, $\ell_0(\vec{x})$ and $\beta(\vec{x})$ are positive, continuous, bounded away from 0 and $\infty$, $\ell_1$ is continuous and $|\ell_1(\vec{x})|$ is bounded away from $\infty$, and $\nu(y|\vec{x})$ is the remaining term satisfying 
\[
\sup_{\vec{x}\in{\cal X}}y^{\beta(\vec{x})}|\nu(y|\vec{x})|\rightarrow 0
 \ {\rm and}\ \sup_{\vec{x}\in{\cal X}} y\left|\frac{\partial \nu(y|\vec{x})}{\partial y}\right| \rightarrow 0
,\ \ {\rm as}\ y\rightarrow\infty.
\]
Because
\[
\frac{\partial L(y|\vec{x})}{\partial y} = -\ell_1(\vec{x})\beta(\vec{x})y^{-\beta(\vec{x})-1} +\frac{\partial \nu(y|\vec{x})}{\partial y},
\]
the density function of (\ref{paretotail}), $f(y|\vec{x})=\partial F(y|\vec{x})/\partial y$, becomes
\begin{eqnarray}
f(y|\vec{x}) &=& \frac{1}{\gamma^*(\vec{x})}y^{-1/\gamma^*(\vec{x})-1}L(y|\vec{x})+ y^{-1/\gamma^*(\vec{x})}\frac{\partial L(y|\vec{x})}{\partial y}\nonumber\\
&=&  \frac{\ell_0(\vec{x})}{\gamma^*(\vec{x})}y^{-1/\gamma^*(\vec{x})-1}\{1+o(1)\} \label{truedensity}
\end{eqnarray}
as $y\rightarrow\infty$. 

As described in this paper, EVI is assumed to be expressible as the single-index model: $\gamma^*(\vec{x})=\gamma(\vec{x}^\top\vec{\theta})$, where $\gamma:\mathbb{R}\rightarrow \mathbb{R}_+$ is the univariate nonlinear function and $\vec{\theta}=(\theta_1,\ldots,\theta_p)^\top\in \mathbb{R}^p$ is the single-index parameter vector. 
However, the pair of the true structure $(\gamma, \vec{\theta})$ is well known not to be unique (Ichimura 1993, Kuchibhotla and Patra 2016). 
To identify this point, we assume that $\vec{\theta}\in {\cal S}_+^{p-1}$, where $\|\cdot\|$ is the Euclidean norm, and 
\[
 {\cal S}_+^{p-1}\equiv \left\{(\theta_1,\ldots,\theta_p)^\top\left| \|\vec{\theta}\|=1, \theta_1\geq 0\right.\right\}.
\]
This assumption identifies the scale and sign of the single index parameter vector. 
It is noteworthy that the constraint of single index parameter vector above excludes the case in which no covariate $\vec{X}$ is predictive for EVI. 
However, the no-covariate model can be characterized by the case in which the nonlinear function $\gamma(\cdot)$ is reduced to the constant in $\vec{x}^\top\vec{\theta}$ for all $\vec{x}\in{\cal X}$ and $\vec{\theta}\in{\cal S}_+^{p-1}$ (as described in Remark 2 of Section 2.2).
Thus, the Pareto-type tailed distribution with the single-index model is defined as 
\begin{eqnarray}
P(Y>y|\vec{X}=\vec{x})= y^{-1/\gamma(\vec{x}^\top \vec{\theta})}L(y|\vec{x}). \label{paretotail2}
\end{eqnarray}

Letting $\{(Y_i,\vec{X}_i):i=1,\ldots,n\}, \vec{X}_i=(X_{i1},\ldots,X_{ip})^\top$ be an $i.i.d.$ random sample generated from a distribution similar to $(Y,\vec{X})$, then to estimate $(\gamma, \vec{\theta})$, we use the POT method.
Subsequently, we introduce threshold $w$ and estimate $(\gamma,\vec{\theta})$ using all observations that exceed the threshold: $\{(Y_i,\vec{X}_i):Y_i>w ,i=1,\ldots,n\}$. 
Actually, given $Y>w$ and $\vec{X}=\vec{x}$, the transformed random variable $Y/w$ is distributed as 
\begin{eqnarray}
P\left(\left.\frac{Y}{w}>z\right|\vec{X}=\vec{x}, Y>w\right)=\frac{1-F(zw|\vec{X}=\vec{x})}{1-F(w|\vec{X}=\vec{x})}= z^{-1/\gamma(\vec{x}^\top\vec{\theta})}\{1+o(1)\},\ \ z\geq 1, \label{paretotail-est}
\end{eqnarray}
as $w\rightarrow\infty$. 
Consequently, the Pareto-type tailed distribution can be replaced approximately with an ordinary Pareto distribution.
Hereinafter, we continue the discussion using (\ref{paretotail-est}).
The conditional density function $f_{w}(\cdot|\vec{x})$ of $Y/w$ given $\vec{X}=\vec{x}$ and $Y>w$ is obtained as
\begin{eqnarray*}
f_{w}\left(\left.\frac{y}{w}\right|\vec{x}\right)
\approx
\frac{1}{\gamma(\vec{x}^\top\vec{\theta}|\vec{\theta})} \left(\frac{y}{w}\right)^{-\displaystyle\frac{1}{\gamma(\vec{x}^\top\vec{\theta})}-1}.
\end{eqnarray*}
Using these expressions, we have
\begin{eqnarray}
-\log f_{w}\left(\left.\frac{y}{w}\right|\vec{x}\right)
&\approx&
\left(\frac{1}{\gamma(\vec{x}^\top\vec{\theta})}+1\right)\log\left(\frac{y}{w}\right)-\log \frac{1}{\gamma(\vec{x}^\top\vec{\theta})} \nonumber\\
&=&
(\exp[\alpha(\vec{x}^\top\vec{\theta})]+1)\log\left(\frac{y}{w}\right)-\alpha(\vec{x}^\top\vec{\theta}), \label{density}
\end{eqnarray}
where $\alpha(\cdot)=-\log \gamma(\cdot)$. 
By expressing $\gamma(\cdot)$ as $\exp[-\alpha(\cdot)]$ and by estimating $\alpha$ instead of $\gamma$ directly, the positivity of $\gamma$ can be ensured naturally.
As described in Section 2.2, we estimate $(\alpha, \vec{\theta})$ using penalized maximum likelihood based on the logarithm of the approximated density function above. 
The estimator of the EVI for the point $\vec{x}\in{\cal X}$ is obtained by $\hat{\gamma}(\vec{x}^\top\hat{\vec{\theta}})=\exp[-\hat{\alpha}(\vec{x}^\top\hat{\vec{\theta}})]$, where $\hat{\alpha}$ is the estimator of $\alpha$ and $\hat{\vec{\theta}}$ is the estimator of $\vec{\theta}$.
The proposed estimator can be regarded as a semiparametric version of the linear estimator proposed by Wang and Tsai (2009).

\noindent{\bf Remark 1} 

As described in this paper, the single-index assumption is incorporated only for EVI $\gamma$. 
This assumption can be extended to the conditional distribution as $F(y|\vec{x})=F(y|\vec{x}^\top\vec{\theta})$ or $L(y|\vec{x})=L(y|\vec{x}^\top\vec{\theta})$. 
However, the information of $L$ is not used to estimate the EVI function by POT. 
Therefore, the single-index assumption for $L$ is unimportant. 
For that reason, we use the single index model only for EVI.
In fact, the sufficient condition of the assumption $F(y|\vec{x})=F(y|\vec{x}^\top\vec{\theta})$ is discussed by Zhu et al. (2012), along with results presented by Li (1991) and by Hall and Li (1993). 
In this sense, one might also naturally assume that $F(y|\vec{x})=F(y|\vec{x}^\top\vec{\theta})$.

\subsection{Estimation procedure}

It is now possible to estimate $(\alpha,\vec{\theta})$ from the data $\{(Y_i,\vec{X}_i):i=1,\cdots,n\}$. 
The nonlinear function $\alpha$ is approximated using the spline method. 
As described herein, we assume that for any $\vec{X}\in{\cal X}$ and $\vec{\theta}\in {\cal S}_+^{p-1}$, there exist $a, b$ such that $a\leq \vec{X}^\top\vec{\theta}\leq b$ (see (C1) in Section 3). 
Let ${\cal C}^{d}[a,b]$ be the class of functions with $d$-times continuously differentiable on $[a,b]$. 
We then define the set of knots $\vec{\kappa}=\{a=\kappa_0<\kappa_1<\ldots<\kappa_{K_0+1}=b\}, K_0>1$ and the class of $d$-th order spline as
\[
{\cal S}(d, \vec{\kappa})=\{s\in{\cal C}^{d-2}[a,b]: s\ is\ a\ polynomial\ of\ degee\ (d-1)\ on\ each\ subinterval\ [\kappa_j,\kappa_{j+1}]\}\ \ d\geq 2.
\]
For $d=1$, ${\cal S}(d, \vec{\kappa})$ is the set of step functions with jumps at each knot. 
When $d=4$, it corresponds to the cubic $B$-spline, which is mainly used for data analysis. 
We approximate $\alpha(\cdot)$ by a $d$-th order spline function $s\in{\cal S}(d,\vec{\kappa})$ for some $d>0$ and set of knots $\vec{\kappa}$. 
For $x\in[a,b]$, let $\vec{B}^{[d]}(x)=(B^{[d]}_1(x),\ldots,B^{[d]}_{K}(x))^\top$ be the vector of $d$-th order scaled $B$-spline basis with $K=K_0+d$ (Appendix A). 
For simplicity, we write $\vec{B}(x)=\vec{B}^{[d]}(x)$ and $B_j(x)=B_j^{[d]}(x)$.
From de Boor (2001), all $d$-th order spline functions can be expressed as linear combinations of $B$-spline bases. 
In other words, for any $s\in{\cal S}(d, \vec{\kappa})$, there exists $\vec{b}=(b_1,\ldots,b_{K})\in\mathbb{R}^{K}$ such that for any $z\in[a,b]$, $s(z)=\vec{B}(z)^\top\vec{b}$. 
From this point, for a fixed $\vec{\theta}\in{\cal S}_+^{p-1}$, $\alpha(\vec{x}^\top\vec{\theta})$ is approximated as $\vec{B}(\vec{x}^\top\vec{\theta})^\top\vec{b}$. 
Let 
\begin{eqnarray}
\ell_n(\vec{b}, \vec{\theta}|\lambda)&=&\frac{1}{n}\sum_{i=1}^n \left[\exp[\vec{B}(\vec{X}_i^\top\vec{\theta})^\top\vec{b}]\log\left(\frac{Y_i}{w}\right)-\vec{B}(\vec{X}_i^\top\vec{\theta})^\top\vec{b}\right]I\left(Y_i>w_n\right) \nonumber\\
&&+\frac{\lambda}{2} \int_a^b \left\{\frac{d^m}{dx^m} \vec{B}(x)^\top\vec{b}\right\}^2 dx \label{likelihood}
\end{eqnarray}
be the penalized (minus) log-likelihood loss function obtained from (\ref{density}), where $\lambda>0$ is the smoothing parameter. 
The estimator of $(\vec{b}, \vec{\theta})$ is defined as 
\[
(\hat{\vec{b}}, \hat{\vec{\theta}}) = \underset{\vec{b}\in\mathbb{R}^{K}, \vec{\theta}\in{\cal S}^{p-1}_+}{\argmin}\ \ \ell_n(\vec{b}, \vec{\theta}|\lambda).
\]
For any $\vec{x}\in{\cal X}$, the EVI function is estimated as $\hat{\alpha}(\vec{x}^\top\hat{\vec{\theta}})= \vec{B}(\vec{x}^\top\hat{\vec{\theta}})^\top\hat{\vec{b}}$. 

In practice, the single index parameter vector should be estimated in ${\cal S}_+^{p-1}$, which engenders difficult optimization.
To avoid such difficulties, we reparameterize $\vec{\theta}$. 
Let ${\cal S}_*=\{(\phi_1,\ldots,\phi_{p-1})\in\mathbb{R}^{p-1}: \|\vec{\phi}\| \leq 1\}$. 
We then write $\vec{\theta}=\vec{\theta}(\vec{\phi})=(\sqrt{1-\|\vec{\phi}\|^2},\vec{\phi}^\top)^\top$ for $\vec{\phi}\in{\cal S}_*$. 
Such a transformation, provided by Yu and Ruppert (2002), describes the condition $\|\vec{\theta}\|=1$ and $\theta_1\geq 0$ with only restriction $\|\vec{\phi}\leq 1$. 
As an alternative to (\ref{likelihood}), we construct the estimator as 
\begin{eqnarray}
(\hat{\vec{b}}, \hat{\vec{\phi}}) = \underset{\vec{b}\in\mathbb{R}^{K}, \vec{\phi}\in{\cal S}_*}{\argmin}\ \ \ell_n(\vec{b}, \vec{\theta}(\vec{\phi})|\lambda) \label{ObjectiveFunc}
\end{eqnarray}
and $\hat{\vec{\theta}}=\vec{\theta}(\hat{\vec{\phi}})$.

\noindent{\bf Remark 2} 

Because the single-index parameter vector is constrained by $\|\vec{\theta}\| = 1$, it might appear that the proposed model inherently excludes the null model (i.e., a model without covariate effects).
However, the null model can still be represented via the nonlinear component.
Specifically, because of the properties of the scaled $B$-spline basis (Appendix A), if $b_k=c K^{-1/2}$ for some constant $c$, then $\vec{B}(z)^\top \vec{b} = c$ for any $z \in [a,b]$.
Regarding estimation from (\ref{likelihood}), when $\lambda \to \infty$, the estimated function $\hat{\alpha}(\cdot)$ is reduced to a polynomial of degree $m-1$ in $\vec{X}^\top \vec{\theta}$.
When the true model is indeed the null model, the estimated slope components in this polynomial are expected to shrink toward zero.

\subsection{Implementation}

Parameters $(\hat{\vec{b}}, \hat{\vec{\phi}})$ are estimated by alternate optimization. 
Let $\vec{\phi}^{(0)}$ be the initial estimator of $\vec{\phi}$. 
At the $k$-th iteration, the updates are given as
\begin{eqnarray*}
\hat{\vec{b}}^{(k)} = \underset{\vec{b}}{\argmin}\ \ \ell(\vec{b},\vec{\phi}^{(k-1)}),
\end{eqnarray*}
and 
\begin{eqnarray}
\hat{\vec{\phi}}^{(k)} = \underset{\|\vec{\phi}\|\leq 1}{\argmin}\ \ \ell(\vec{b}^{(k)},\vec{\phi}). \label{NRphi}
\end{eqnarray}
The iteration continues until $\|\vec{\phi}^{(k)}-\vec{\phi}^{(k-1)}\|<\varepsilon$ for some $\varepsilon>0$. 
As described herein, we set $\varepsilon=10^{-4}$. 
At each step, $\hat{\vec{b}}^{(k)}$ is computed using the {\sf optim} function in R, whereas the optimization of $\hat{\vec{\phi}}^{(k)}$ with norm constraint $\|\vec{\phi}^{(k)}\|\leq 1$ is obtained via the {\sf constrOptim.nl} function in the {\sf alabama} package (see Varadhan 2023). 
To accelerate the estimation of $\vec{\phi}$ further, we suggest use of the proximal descent algorithm, modifying (\ref{NRphi}) to
\begin{eqnarray*}
\hat{\vec{\phi}}^{(k)} = \underset{\vec{\phi}}{\argmin}\ \ \ell(\vec{b}^{(k)},\vec{\phi}) + \nu^{(k)}\|\vec{\phi}-\vec{\phi}^{(k-1)}\|^2,
\end{eqnarray*}
where $\nu^{(k)}$ represents the step size. 
In our numerical experiments, setting $\nu^{(0)}=10^{-5}$ and updating as $\nu^{(k)}=2 \nu^{(k-1)}$ yielded fast and stable convergence. 

The algorithm in this section is implemented under a fixed tuning parameter setting $(w_n,\lambda)$.
In practice, the final estimator is obtained by selecting tuning parameters as described in Section 4.2.
In our experiments, the initial value $\vec{\phi}^{(0)}$ was chosen via multiple random starts for the first tuning parameter configuration $(w_n,\lambda)$. 
For subsequent configurations, we adopted a warm-start strategy using the estimate obtained from the previous tuning parameter as the initial value.
It is noteworthy that the unit vector cannot be used as the initial $\vec{\phi}^{(0)}$ because it implies the boundary of the condition $\|\vec{\phi}^{(0)}\|\leq 1$.

\subsection{Tuning parameter selection}

In the proposed estimator, we have the following three tuning parameters as the threshold $w_n$, number of knots $K$, and smoothing parameter $\lambda$. 
According to Ruppert (2002), knots selection is not as important as $\lambda$. 
Ruppert demonstrated that using equidistant knots with fixed large $K$ is sufficient.
For our method, we choose $w_n$ and $\lambda$ using the data-driven method. 
To choose $w_n$, we use the discrepancy measure provided by Wang and Tsai (2009). 
Letting $U_i=\exp[-\exp[\alpha(\vec{X}_i^\top\vec{\theta})]\log(Y_i/w_n)]$, then $U_i$ approximately distributed to a standard uniform distribution under $Y_i>w_n$. 
Therefore, the criterion of goodness of fit can be used to the standard uniform distribution to detect the tuning parameter. 
Actually, we use $\hat{U}_i= \exp[-\exp[\hat{\alpha}(\vec{X}_i^\top\hat{\vec{\theta}})]\log(Y_i/w_n)]$.
Define $n_0= \sum_{i=1}^n I(Y_i>w_n)$.
We then define the discrepancy measure as  
\[
D(w_n|\lambda)=\frac{1}{n_0}\sum_{i=1}^{n_0} \{\hat{U}_{(i)}-\hat{F}(i/(n_0+1))\}^2,
\]
where $\hat{U}_{(1)}\leq\cdots\leq \hat{U}_{(n_0)}$ are order statistics of $\{\hat{U}_1,\ldots,\hat{U}_{n_0}\}$ and $\hat{F}(u)$ is the empirical distribution based on $\{\hat{U}_1,\ldots,\hat{U}_{n_0}\}$. 
Tuning parameters $w_n$ are selected via minimizing $D(w_n|\lambda)$ given $\lambda$. 

Next, $\lambda$ is selected by $H$-fold cross validation. 
Dataset $\{(Y_i, \vec{X}_i): i=1,\ldots,n\}$ is partitioned randomly into $H$ disjoint subsets ${\cal J}_1,\ldots,{\cal J}_H$. 
Let $(\hat{\vec{b}}^{[-h]},\hat{\vec{\phi}}^{[-h]})$ be the estimator obtained by (\ref{ObjectiveFunc}) using data excluding those data within ${\cal J}_h$. 
Then, the evaluation score of cross-validation is defined as 
\[
{\rm CV}(\lambda|w_n) = \frac{1}{H}\sum_{h=1}^H \frac{1}{|{\cal J}_h|}\sum_{(Y_i, \vec{X}_i)\in{\cal J}_h} \ell_i(\hat{\vec{b}}^{[-h]}, \hat{\vec{\phi}}^{[-h]}|w_n),
\]
where 
\[
\ell_i(\vec{b},\vec{\phi}|w)=  \left\{\exp[\vec{B}(\vec{X}_i^\top\vec{\theta}(\vec{\phi}))^\top\vec{b}]\log\left(\frac{Y_i}{w}\right)-\vec{B}(\vec{X}_i^\top\vec{\theta}(\vec{\phi}))^\top\vec{b}\right\}I(Y_i>w).
\]
In our numerical study of Sections 4.3 and 5, we used $H=5$.

The selecting algorithm proceeds as follows. Letting $\{w_1,\ldots,w_S\}$ be the set of candidate threshold values and letting$\{\lambda_1,\ldots,\lambda_T\}$ be the set of candidate smoothing parameters, then for each $s=1,\ldots,S$, we calculate $\lambda_{s,cv}=\argmin_t {\rm CV}(\lambda_t|w_s)$. 
Subsequently, the optimal tuning parameters are selected as $(w_{s^*}, \lambda_{s^*,cv} )$ with $s^*=\argmin_s D(w_s|\lambda_{s,cv})$.

\section{Asymptotic Theory}

As described in this section, we study the asymptotic property of the proposed estimator. 
The true parameter and function in (\ref{paretotail2}) is defined as $\vec{\theta}_0$ and $\gamma_0$.
That is, 
\begin{eqnarray*}
P(Y>y|\vec{X}=\vec{x})= y^{-1/\gamma_0(\vec{x}^\top \vec{\theta}_0)}L(y|\vec{x}).
\end{eqnarray*}
Additionally, we define $\alpha_0(\cdot)=-\log \gamma_0(\cdot)$ and $\vec{\phi}_0$ as $\vec{\theta}_0=\vec{\theta}(\vec{\phi}_0)$.

We consider the following conditions.
\begin{enumerate}
\item[(C1)] The marginal density function of $\vec{X}$ is continuous and bounded away from 0 and $\infty$. 
The support ${\cal X}$ of $\vec{X}$ is compact. 
Furthermore, there exist $a,b\in\mathbb{R}$ such that for any $\vec{X}\in{\cal X}$ and any $\vec{\theta}\in{\cal S}_+^{p-1}$, $a\leq \vec{X}^\top\vec{\theta}\leq b$. 

\item[(C2)] $\alpha_0\in {\cal C}^{q}[a,b]$ for some positive $q$. For the order of spline $d$ and the order of the difference penalty $m$ in (\ref{likelihood}), $m < d\leq q$. 

\item[(C3)] In (\ref{slowly}), constant $\beta_{inf}>0$ exists such that $\beta_{inf}\leq \inf_{\vec{x}\in{\cal X}} \gamma(\vec{x}^\top\vec{\theta}_0)\beta(\vec{x})$.

\item[(C4)] The threshold value $w=w_n$ takes $w_n\rightarrow \infty$ such that $\tau_n=E[P(Y>w_n|\vec{X})]$ satisfies $\tau_n\rightarrow 0$ and $n\tau_n\rightarrow \infty$ as $n\rightarrow\infty$.

\item[(C5)] The knots sequence $\vec{\kappa}$ is quasi-uniform: $c_\ell<\max_j\{\kappa_{j+1}-\kappa_j\}/\min_j\{\kappa_{j+1}-\kappa_j\}< c_u $ for some constant $c_\ell, c_u>0$. The number of knots satisfies $K\rightarrow\infty$, but $K\{\log n\}^2/(n\tau_n )\rightarrow 0$ and $K^d\tau_n^{-\beta_{inf}}\rightarrow\infty$ as $n\rightarrow\infty$.

\item[(C6)] The smoothing parameter $\lambda=\lambda_n$ satisfies $\lambda\rightarrow 0$, $\lambda/\tau_n\rightarrow 0$, and $K(\lambda/\tau_n)^{(1/2m)}=O(1)$ as $n\rightarrow\infty$. 
\end{enumerate}

Condition (C1) is a natural condition in nonparametric or semiparametric regression (e.g. Tsybakov 2009). 
For datapoints $\vec{x}_i (i=1,\ldots,n)$ and any $\vec{\theta}\in{\cal S}_+^{p-1}$, we have $-\|\vec{x}_i\|\leq \vec{x}_i^\top\vec{\theta}\leq \|\vec{x}_i\|$. 
Consequently, for example, by centering and scaling the data, $a$ and $b$ are identifiable in practice. 
Wang and Yang (2009) and Wang and Tsai (2009) proposed another method to transform $\vec{X}$ which has known finite support. 
Actually, (C2) is common in the spline smoothing (Xiao 2019). 
In (\ref{likelihood}), the penalty is added to the $m$-th derivative of $\alpha_0$. 
Since $d$th order spline is the $(d-1)$th piecewise polynomial, $m\leq d-1$ and $d\leq q$ are natural. 
In (C3), the value $\gamma(\vec{x}^\top\vec{\theta}_0)\beta(\vec{x})$ can be regarded as a second order parameter in extreme value theory (Section 2, de Haan and Ferreira 2006). 
Together with the Hall class assumption, the positivity of the second-order parameter $\beta_{inf}$ is natural. 
Next we consider (C4). 
The number of data exceeding the threshold is $n_0=\sum_{i=1}^n I(Y_i>w_n)$; we obtain $E[n_0]/n= P(Y_i>w_n)=E[P(Y>w_n|\vec{X})]=\tau_n$.
Consequently, $\tau_n$ controls the rate of data exceeding the threshold. Also, $n\tau_n$ can be regarded as the effective sample size. 
Condition (C4) means that the effective sample size becomes large, but its rate is lower than the original sample size $n$. 
This rate of the effective sample size is called the intermediate order sequence in extreme value theory (Section 2, de Haan and Ferreira 2006).
Actually, (C5) is necessary to obtain a good $B$-spline estimator of the true nonlinear function. This is the standard setting for the $B$-spline method (Xiao 2019). 
The condition $K=o(n\tau_n/\{\log n\}^2)$ indicates that the number of knots cannot be greater than the effective sample size. 
The term $(\log n)^2$ is necessary to prove Lemma \ref{Consistency} rigorously.
Roughly speaking, $O(K^{-d})$ is the rate of approximation bias of the spline function (Lemma \ref{appspline} of Appendix B), whereas $O(\tau_n^{\beta_{inf}})$ is the order of bias resulting from approximating the Pareto distribution (\ref{paretotail2}). 
The approximation bias of the Pareto distribution is related to the second-order condition in EVT (Section 2.2, de Haan and Ferreira 2006), which cannot be ignored because the occurrence of such bias is a common problem in EVT. 
The bias of the spline approximation is dominated by the bias from the penalty term in the penalized spline method (Xiao 2019), as reflected in (C6). 
Therefore, we assume the condition $K^d\tau_n^{-\beta_{inf}}\rightarrow\infty$ so that the bias of the spline model approximation is of negligible order compared to that of the Pareto tail distribution. 
Actually, (C6) is important for penalized spline smoothing, which is related to Remark 5.3(b) and Remark 6.6 reported by Xiao (2019).
If (C6) is violated, then the estimator might not be a consistent estimator of the true nonlinear function.

We can let 
\[
\vec{b}_0= \underset{\vec{b}\in \mathbb{R}^K}{\argmin} \ L(\vec{b}),
\]
where
\[
L(\vec{b})= E\left[\exp[\vec{B}(\vec{X}^\top \vec{\theta}_0)^\top\vec{b}]\log(Y/w_n)-\vec{B}(\vec{X}^\top \vec{\theta}_0)^\top\vec{b}|Y>w_n\right].
\]
Lemma 1 in Appendix B shows that $\sup_{x\in[a,b]}|\alpha_0(x)-\vec{B}(x)^\top\vec{b}_0|=O(K^{-d
})$, which is the optimal asymptotic rate of the spline approximation. 
In other words, $\vec{b}_0$ is the coefficient of best approximation of $B$-spline function to $\alpha_0$. 
One can recall that $\vec{\phi}_0$ satisfies $\vec{\theta}_0=\vec{\theta}(\vec{\phi}_0)$.
First, we show the asymptotic rate of the estimators $\hat{\vec{b}}$ and  $\hat{\vec{\phi}}$. 

\begin{theorem}\label{RateParameter}
Presuming that (C1)--(C6), then
as $n\rightarrow\infty$,
\begin{eqnarray*}
E[\|\hat{\vec{b}}-\vec{b}_0\|^2]\leq O((n\tau_n)^{-1}(\lambda/\tau_n)^{-1/(2m)}) + O(\lambda/\tau_n)+ O(\tau_n^{2\beta_{inf}}).
\end{eqnarray*}
Under the condition $\lambda/\tau_n= (n\tau_n)^{-2m/(2m+1)}$, 
\[
E[\|\hat{\vec{b}}-\vec{b}_0\|^2] \leq O((n\tau_n)^{-2m/(2m+1)})+ O(\tau_n^{2\beta_{inf}}).
\]
For the part of single-index parameter vector, as $n\rightarrow\infty$,
\begin{eqnarray*}
E[\|\hat{\vec{\phi}}-\vec{\phi}_0\|^2] \leq O((n\tau_n)^{-1})+ O(\tau_n^{2\beta_{inf}}).
\end{eqnarray*}
\end{theorem}

The first result of theorem \ref{RateParameter} implies that 
\[
\int_a^b \{\hat{\alpha}(z)-\alpha_0(z)\}^2dz \leq O((n\tau_n)^{-2m/(2m+1)})+ O(\tau_n^{2\beta_{inf}})
\]
under $\lambda/\tau_n= (n\tau_n)^{-2m/(2m+1)}$.
Consequently, although it is not surprising, the rate of convergence of the parametric part is faster than that of the semiparametric part. 
The term $O((n\tau_n)^{-2m/(2m+1)})$ can be regarded as the optimal rate of the semiparametric estimator with sample size $n\tau_n$ (e.g., Tsybakov 2009). 
This observation implies that the asymptotic convergence rate of the single-index estimator of the EVI function is dominated by the semiparametric inference, as stated in the theorem below. 

\begin{theorem}\label{RateEstimator}
Presuming (C1)--(C6), then, as $n\rightarrow\infty$,
\begin{eqnarray*}
E\left[\left\{\hat{\alpha}(\vec{X}^\top \hat{\vec{\theta}})- \alpha_0(\vec{X}^\top\vec{\theta}_0)\right\}^2 \right]
\leq  O((n\tau_n)^{-1}(\lambda/\tau_n)^{-1/(2m)}) + O(\lambda/\tau_n)+ O(\tau_n^{2\beta_{inf}}).
\end{eqnarray*}
Under the condition $\lambda/\tau_n= (n\tau_n)^{-2m/(2m+1)}$, 
\begin{eqnarray*}
E\left[\left\{\hat{\alpha}(\vec{X}^\top \hat{\vec{\theta}})- \alpha_0(\vec{X}^\top\vec{\theta}_0)\right\}^2 \right]
\leq O((n\tau_n)^{-2m/(2m+1)})+ O(\tau_n^{2\beta_{inf}}). 
\end{eqnarray*}
In addition, if $\tau_n$ can be taken as $O(n^{-\{2m/(2m+1)\}/\{2m/(2m+1)+2\beta_{inf}\}})$, then the optimal rate of convergence is
\begin{eqnarray*}
E\left[\left\{\hat{\alpha}(\vec{X}^\top \hat{\vec{\theta}})- \alpha_0(\vec{X}^\top\vec{\theta}_0)\right\}^2 \right]
\leq O(n^{-\frac{2\beta_{inf}}{2\beta_{inf}+1 + 1/m}}).
\end{eqnarray*}
\end{theorem}

The rate of convergence in Theorem \ref{RateEstimator} is independent of the dimension of covariate $p$, which indicates that the curse of dimensionality can be avoided. 
For comparison, Goegebeuer et al. (2015) developed the rate of convergence of the fully nonparametric estimator of the EVI function, but its rate becomes lower with the dimension of covariates.

Next, we express $\tau_n=k/n$ with some $k$ satisfying $k\rightarrow \infty$ and $k/n\rightarrow 0$. 
We also write $\rho=-\beta_{inf}$. 
Then, in the theorem \ref{RateEstimator}, the second assertion is 
\begin{eqnarray}
E\left[\left\{\hat{\alpha}(\vec{X}^\top \hat{\vec{\theta}})- \alpha_0(\vec{X}^\top\vec{\theta}_0)\right\}^2 \right]
\leq O(k^{-2m/(2m+1)})+ O((n/k)^{2\rho}) \label{Optimal1} 
\end{eqnarray}
and the last assertion becomes 
\begin{eqnarray}
E\left[\left\{\hat{\alpha}(\vec{X}^\top \hat{\vec{\theta}})- \alpha_0(\vec{X}^\top\vec{\theta}_0)\right\}^2 \right]
\leq O(n^{\frac{2\rho}{1-2\rho + 1/m}}). \label{Optimal2}
\end{eqnarray}

One might also consider the Pareto-type distribution with the Hall class of one-dimensional data as $Y_1,\ldots,Y_n\stackrel{i.i.d.}{\sim} F$, where 
\[
P(Y>y)=1-F(y)=y^{-1/\gamma}\{\ell_0+\ell_1y^{-\beta}+\nu(y)\}
\]
with $\gamma, \ell_0, \beta>0$, $\ell_1\in\mathbb{R}$, $\nu(y)=o(y^{-\beta})$ and $y \partial \nu(y)/\partial y =o(1)$ as $y\rightarrow\infty$. 
Then, the maximum likelihood estimator of $\gamma$ is 
\[
\hat{\gamma}_{{\rm Hill}}=\frac{1}{n}\sum_{i=1}^{n}\log\left(\frac{Y_{i}}{w_n}\right)I(Y_i>w_n),
\]
which is fundamentally similar to the Hill estimator (Hill 1975).
From de Haan and Ferreira (2006), the Hill estimator $\hat{\gamma}_{{\rm Hill}}$ is well known to have an asymptotic rate of convergence as 
\begin{eqnarray}
E[|\hat{\gamma}_{{\rm Hill}}-\gamma|^2]=O(k^{-1}) + O((n/k)^{2\rho}) \label{Hillone}
\end{eqnarray}
under some suitable conditions, $\rho=-\gamma\beta$, $k=\sum_{i=1}^n I(Y_i>w_n)$, $nP(Y>w_n)\rightarrow\infty$ and $P(Y>w_n)\rightarrow 0$ as $n\rightarrow\infty$. 
Parameter $\rho$ is the second-order parameter (Section 2, de Haan and Ferreira 2006). 

Compared to (\ref{Optimal1}) and (\ref{Hillone}), the first terms of both estimators represent the difference between the optimal asymptotic order of the semiparametric estimator and the parametric estimator with sample size $k$.
The second terms are inherently similar. 
Drees (2001) presents the optimal rate of convergence of $\hat{\gamma}_{{\rm HILL}}$ as 
\[
E[|\hat{\gamma}_{{\rm Hill}}-\gamma|^2]=O\left(n^{\frac{2\rho}{1-2\rho}}\right),
\]
which is a slightly higher rate than that reported from an earlier study (\ref{Optimal2}). 
This finding indicates that the results of Theorem \ref{RateEstimator} are a natural extension from the parametric method to the semiparametric regression.

\section{Numerical Experiments}

In this section, we investigate the finite-sample performance of the proposed estimator through Monte Carlo simulations.
For response $Y$ and the covariate $\vec{X}$, the true distribution is set as 
\begin{eqnarray}
P(Y>y|\vec{X}=\vec{x})= \frac{y^{-1/\gamma^*(\vec{x})}}{1+\ell y^{-1/\gamma^*(\vec{x})} }, \label{True}
\end{eqnarray}
where $\ell$ is the positive constant and $\gamma^*:\mathbb{R}^p\rightarrow\mathbb{R}$ is the EVI function. 
The model (\ref{True}) is obtained by (\ref{paretotail}) with $L(y|\vec{x})=(1+\ell y^{-1/\gamma^*(\vec{x})} )^{-1}$. 
It is readily apparent that above $L(y|\vec{x})$ has the form (\ref{slowly}) by setting $\ell_0(\vec{x})=1$, $\ell_1(\vec{x})=-\ell$, $\beta(\vec{x})=1$ and $\nu(y|\vec{x})=\ell^2y^{-2/\gamma^*(\vec{x})}(1+o(1))$ as $y\rightarrow\infty$. 
Because 
\[
\frac{\partial L(y|\vec{x})}{\partial y}= \frac{-\ell}{\gamma^*(\vec{x})}y^{-1/\gamma^*(\vec{x})-1}L(y|\vec{x})^2,
\]
the density function from (\ref{True}) can be written as
\begin{eqnarray}
f(y|\vec{x})&=&\frac{1}{\gamma^*(\vec{x})}y^{-1/\gamma^*(\vec{x})-1} L(y|\vec{x})^2\nonumber\\
&=& 
\frac{1}{\gamma^*(\vec{x})}y^{-1/\gamma^*(\vec{x})-1} (1-\ell y^{-1/\gamma^*(\vec{x})} + O(\ell y^{-2/\gamma^*(\vec{x})} ) ),\ \ {\rm as}\ \ y\rightarrow\infty, \label{trueApp}
\end{eqnarray}
which corresponds to (\ref{truedensity}). 

For our simulation, the covariate $\vec{X}_i=(X_{i1},\ldots,X_{ip})^\top$ is generated as presented below. 
First, for $\vec{Z}_i=(Z_{i1},\ldots,Z_{ip})$, we generate $\vec{Z}_1,\ldots,\vec{Z}_n\sim N(\vec{0},\Sigma)$, where $\Sigma=(0.25^{|k-j|})_{kj}$ for $k,j=1,\ldots,p$. 
Next, we construct $X_{ij}=(1/\sqrt{3})\{\hat{F}_{Z,j}(Z_{i,j-1})-1/2\}$ for $j=1,\ldots,p$, where $\hat{F}_{Z,j}$ denotes the empirical distribution function based on $\{Z_{1j},\ldots,Z_{nj}\}$ for $j=1,\ldots,p$.  
Roughly speaking, for each $j$, $X_{ij}$ is distributed approximately as a uniform distribution on an interval $[-1/\sqrt{3}, 1/\sqrt{3}]$, which implies that $X_{ij}$ has mean 0 and variance 1. 
Furthermore, each pair $(X_{ij},X_{ik})$ can be found to have some correlation. 
Under given $\vec{X}=\vec{x}$, the response $Y$ is generated from (\ref{True}) by inversion. 
Hereinafter, this report describes results obtained from the simulation under some settings.

\subsection{Effect of high-dimensionality of covariates}

This section presents an illustration of the finite sample performance of the estimator varying with the number of covariates.
First, the parameter $\vec{\theta}=(\theta_1,\ldots,\theta_p)^\top\in\mathbb{R}^p$ is prepared, where $\theta_1=1, \theta_2=0.2$ and $\theta_3=0.5$ and $\theta_j=0, j>3$. 
The $\vec{\theta}$ is modified as $\vec{\theta}/\|\vec{\theta}\|$. 
Hereinafter, as the true EVI function $\gamma^*(\vec{x})=\exp[-\alpha^*(\vec{x})]$, we use $\alpha^*(\vec{x})=\alpha(\vec{x}^\top\vec{\theta})=-3+\phi(\vec{x}^\top\vec{\theta};-\mu,\sigma)+ \phi(\vec{x}^\top\vec{\theta};\mu,\sigma)$ with $\mu=0.3$ and $\sigma=0.2$, where $\phi(z;\mu,\sigma)$ is the density function of Gaussian distribution with mean $\mu$ and standard deviation $\sigma$. 
In (\ref{True}), we set $\ell=0.25$ . 
The sample size is fixed as $n=2000$. 
The number of coveriates is set as $p=3, 20$ and 50. 
Consequently, when $p=3$, it corresponds to the true setting. The model with $p>3$ includes the irrelevant covariates.

\begin{figure}
\centering
\includegraphics[width=150mm,height=80mm]{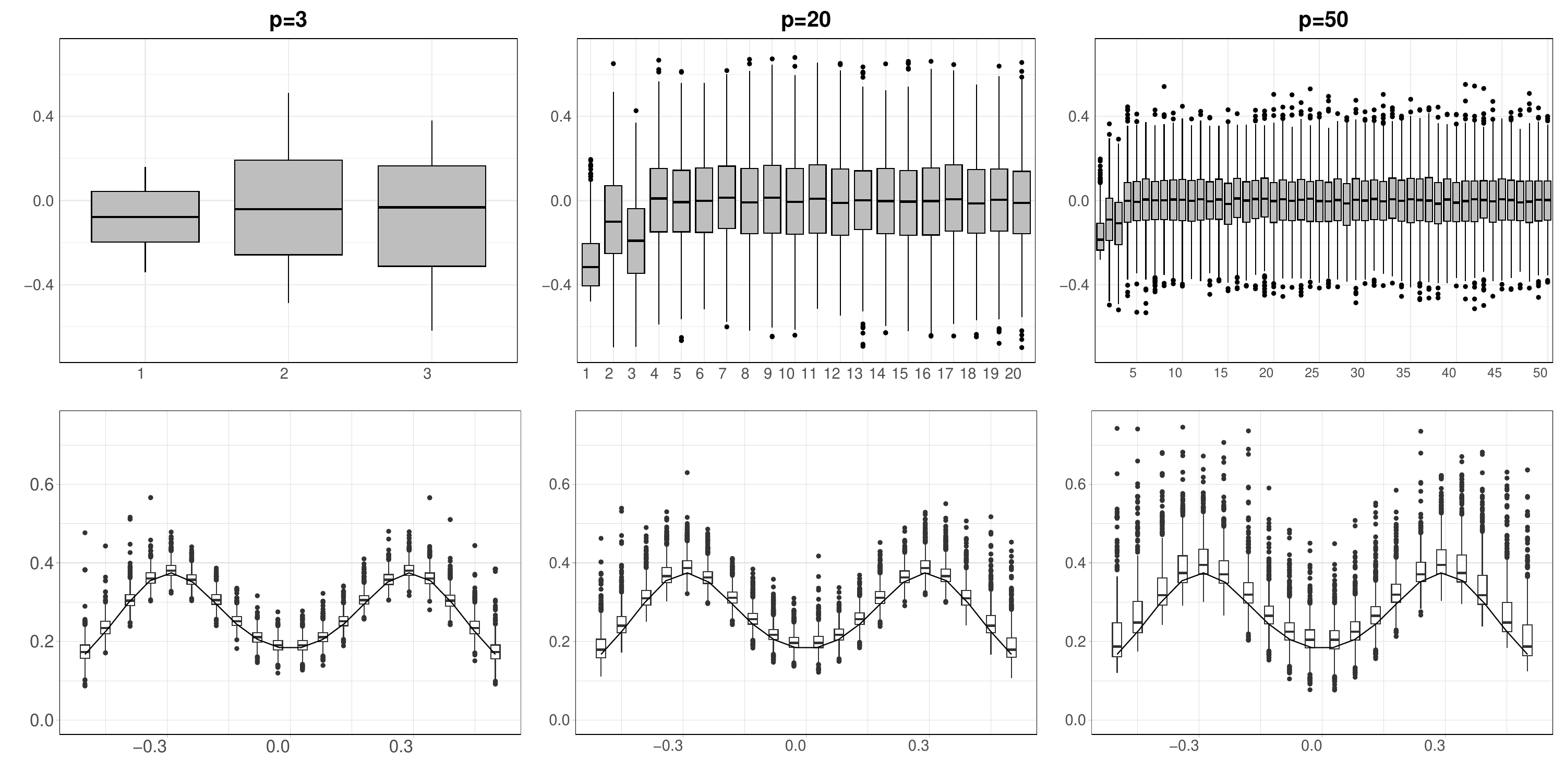}
\centering
\caption{Simulation results corresponding to Section 4.3.  
Upper panels show boxplots of each $\hat{\theta}_j - \theta_j$ for $j = 1, \ldots, p$.  
Lower panels show boxplots of $\hat{\gamma}(z)$ at each $z \in [-0.5, 0.5]$ with line $\gamma(z)$.  
From left to right, the results for $p = 3$, 20, and 50 are presented. \label{Result_Section 43} }
\end{figure}

Figure \ref{Result_Section 43} portrays boxplots of the proposed estimator based on 500 Monte Carlo iterations.
For $p = 3$, both the parametric and semiparametric components performed well.
As $p$ increases, the performance of the estimator deteriorates, which is expected because of the higher dimensionality.
Even for $p = 20$, the estimators capture the underlying structure of the true model, although some bias can be observed.
For $p = 50$, both the parametric components ${\hat{\theta}_j - \theta_j}$ ($j = 1, 2, 3$) and the semiparametric estimator exhibit noticeable bias, probably because of the inclusion of several irrelevant covariates.
Nevertheless, the overall structure of the true model is reasonably well captured.
These results suggest that large values of $p$ tend to engender underestimation.
In future work, we intend to develop a dimension reduction technique to mitigate such effects in very high-dimensional settings.

\subsection{Effect of sample size}

\begin{figure}
\begin{center}
\includegraphics[width=150mm,height=80mm]{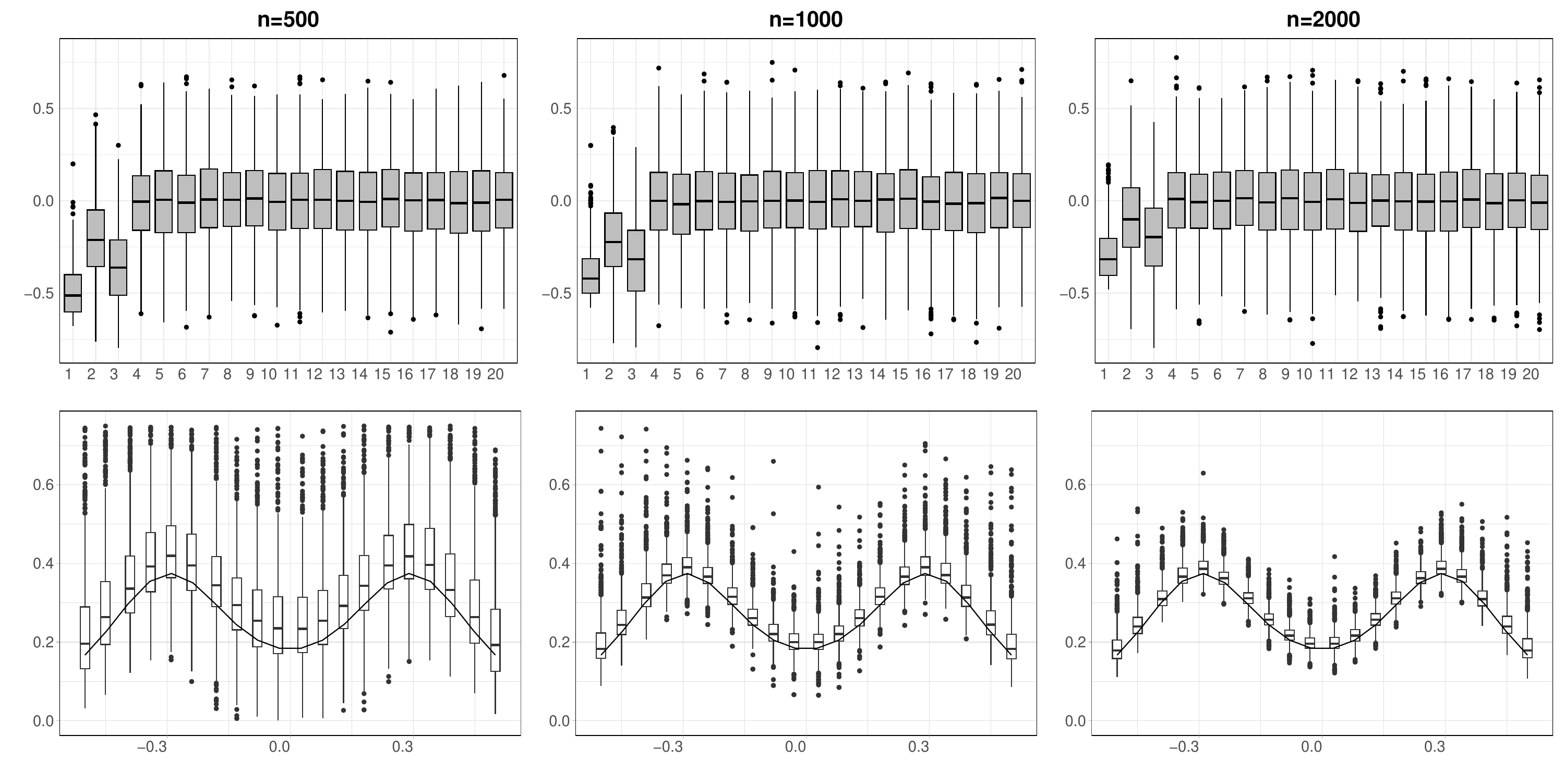}
\end{center}
\caption{Simulation results corresponding to Section 4.4.  
The upper panels present boxplots of each $\hat{\theta}_j - \theta_j$ for $j = 1, \ldots, 20$.  
The lower panels show boxplots of $\hat{\gamma}(z)$ at each $z \in [-0.5, 0.5]$ with line $\gamma(z)$.  
From left to right, results for $n = 500$, 1000, and 2000 are presented.
\label{Result_Section44} } 
\end{figure}

Consider a similar model to that of Section 4.3.2 with $p=20$. 
Here, we confirm the performance of the estimator varying with sample size as $n=500$, $1000$ and $2000$. 
Figure \ref{Result_Section44} presents the boxplot of the performance of the proposed estimator for each $n$. 
For each $n$, there were biases of the estimators of the non-zero components ($\theta_1, \theta_2$ and $\theta_3$), but these biases were reduced by increasing $n$. 
For zero-components $\theta_j, j>3$, the medians of boxplots were close to zero, but the deviances were large, even for $n=2000$. 
For the semiparametric part, the performance of the estimator with $n=500$ was not good because the effective sample size $n_0= \sum_{i=1}^n I(Y_i>w)$ is small. 
Actually, $n_0$ was about average 71.2 with standard error 9.84 for $n=500$ among 500 Monte Carlo replications.  
When $n=1000$ and 2000, it can be said that acceptable results were obtained. 
The average (standard error) of effective sample sizes of $n=1000$ and $2000$ were, respectively, 116.2 (14.35) and 159 (17.46).

\subsection{Robustness to model misspecification}

We apply the proposed method to five models. 
The models are (i) $\alpha^*(\vec{x})=1.2+ 2\vec{x}^\top\vec{\theta}$ and $\ell=0$, (ii) similar to $\alpha^*$ as (i) but $\ell_1=0.25$, (iii) similar model to that in presented Section 4.3.2 and  (iv) $\alpha^*(\vec{x})= -1.2-0.5(1-x_3)\sin(2\pi x_2)$ and $\ell=0.25$. 
The sample size and the number of covariates are fixed respectively as $n=2000$ and $p=20$. 
The estimator performance is evaluated by approximated integrated squared error as
\[
{\rm ISE}= \frac{1}{J}\sum_{j=1}^J \left\{\frac{\hat{\gamma}(\vec{X}_j^*)}{\gamma^*(\vec{X}_j^*)}-1\right\}^2 , 
\]
where $\vec{X}_j^*$ are test data generated from the similar distribution of $\vec{X}$, and where $\hat{\gamma}(\vec{x})=\exp[-\hat{\alpha}(\vec{x}^\top\hat{\vec{\theta}})]$ is the estimator of $\gamma^*(\vec{x})$. 
However, we excluded test data for which $(\vec{X}_j^*)^\top\hat{\vec{\theta}}$ falls outside the central 90\% interval, i.e., below the fifth quantile or above the 95th quantile. The total number of test data is adjusted as $J=1000$.
For our study, we evaluate the distribution of ISE calculated using 500 Monte Carlo iterations. 

\begin{figure}
\begin{center}
\includegraphics[width=120mm,height=60mm]{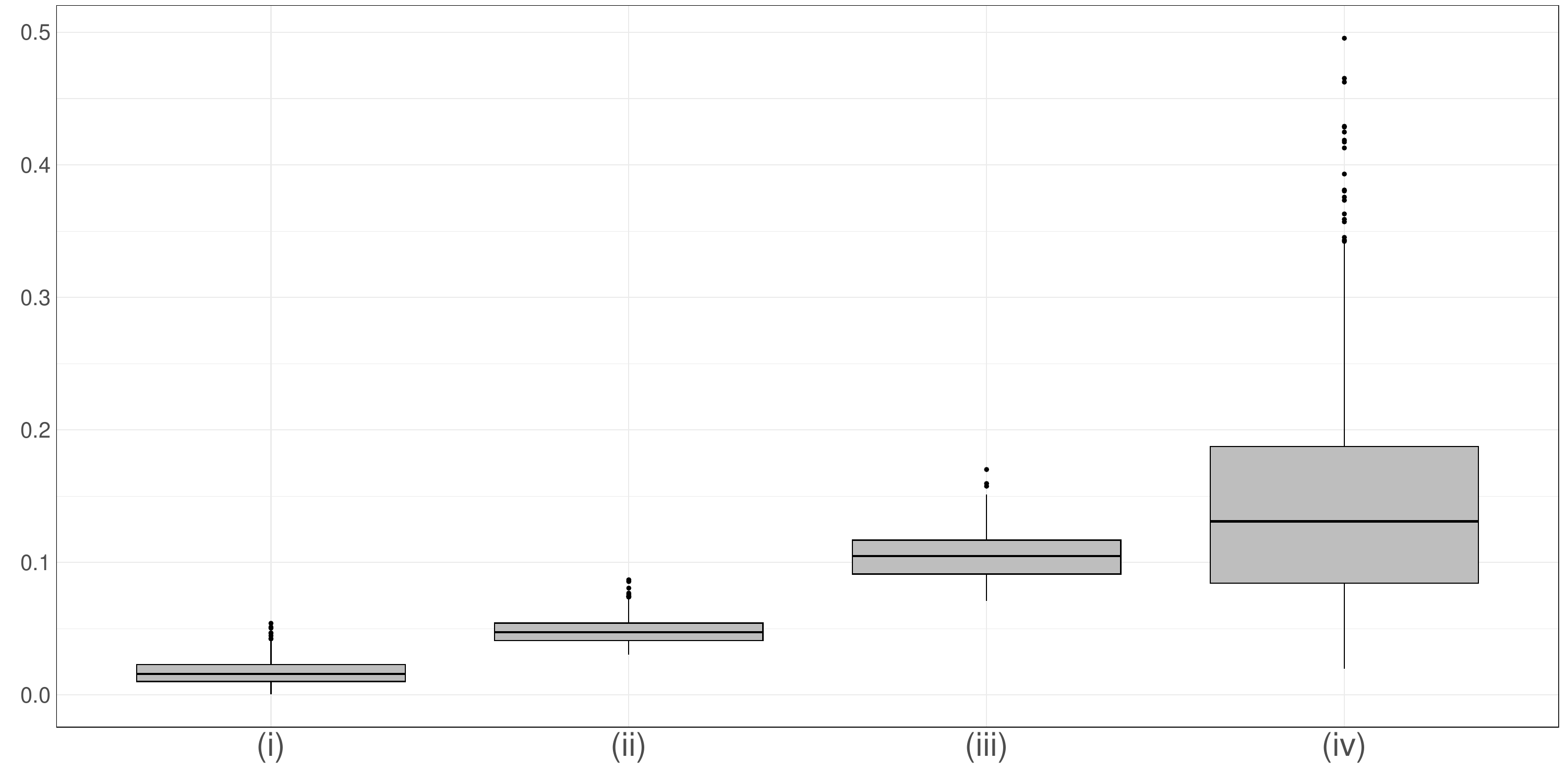}
\end{center}
\caption{Simulation results corresponding to Section 4.5.  
Boxplots of the ISE of the estimator for four models are presented.
\label{Result_Section45} }  
\end{figure}

Results are presented in Figure \ref{Result_Section45}.
Models (i) and (ii) present a quite simple structure providing good performance. 
Models (i) and (ii) have similar EVI function, but the conditional distribution of $Y$ given $\vec{X}=\vec{x}$ differs because of $\ell$ in (\ref{trueApp}). 
Roughly speaking, for (i), the threshold selection is needless but the estimator is constructed with threshold selection. 
Consequently, the effective sample size for (i) tends to be larger than that for (ii). 
From this, ISE for (i) is smaller than that for (ii). 
The model (iii), as detailed in Sections 4.3 and 4.4, has a complicated structure of $\alpha(\cdot)$, but $\ell$ is similar to (ii). 
It is apparent from complexity of $\alpha$ that the estimator performance is worse than that for (ii). 
However, the distribution of the estimator was robust. 
The last model (iv) is a fully nonlinear model with no single-index structure such as those of (i)--(iii). 
It is apparent from the result that the dispersion of the ISE is large. 
However, the median of ISE was approximately equal to that for (iii), which indicates that the single index model is an efficient approach even for the fully nonlinear model.

\subsection{Comparison of various methods} 

This section presents a comparison of the proposed method and other estimator for a model similar to that in Section 4.3 with $n=2000$ and $p=20$. 
As described herein, the proposed estimator with single index model is denoted as SIM. 
The competitors are the following. 
First, the single index model with the tuning parameters $(w,\lambda)$ selected by minimizing ISE, which denotes the Oracle. 
The Oracle is the optimal estimator from our model, but it is calculable only through simulation because the information of true EVI function is used. 
Next, we consider the no-covariate model, which is denoted by Null. 
It is noteworthy that the Null estimator is fundamentally similar to the Hill estimator (Hill 1975).
We use the linear model $\alpha(\vec{x})=\theta_0+\vec{x}^\top\vec{\theta}$, as proposed by Wang and Tsai (2009). 
The additive model $\alpha(\vec{x})=\alpha_0+\alpha_1(x_1)+\cdots+\alpha_p(x_p), \alpha_j:\mathbb{R}\rightarrow\mathbb{R}$ (Youngman 2019) is also considered. 
Then, each $\alpha_j$ is estimated using the smoothing spline method.
All smoothing parameters for $\alpha_1,\ldots,\alpha_p$ are similar. The threshold value for POT selected by the discrepancy measure is described in Section 4.2. 
The estimators using linear and additive models are denoted as Linear and Additive. 
We also considered the method of the quantile-based single index model proposed by Xu et al. (2022), which is denoted by QSIM. 
Let $Q_Y(\tau|\vec{x})$ be the conditional quantile of $Y$, given $\vec{X}=\vec{x}$. 
Their method assumes that $Q_Y(\tau|\vec{x})=Q_Y(\tau|\vec{x}^\top\vec{\theta})=\vec{x}^\top\vec{\theta}(\tau)$.
Then, $\vec{\theta}(\tau)$ is estimated using linear quantile regression.
The EVI is estimated using  the conditional Hill estimator from the estimator of conditional quantile. 
The quantile level is fixed as $\tau=0.9$.

\begin{figure}
\begin{center}
\includegraphics[width=120mm,height=60mm]{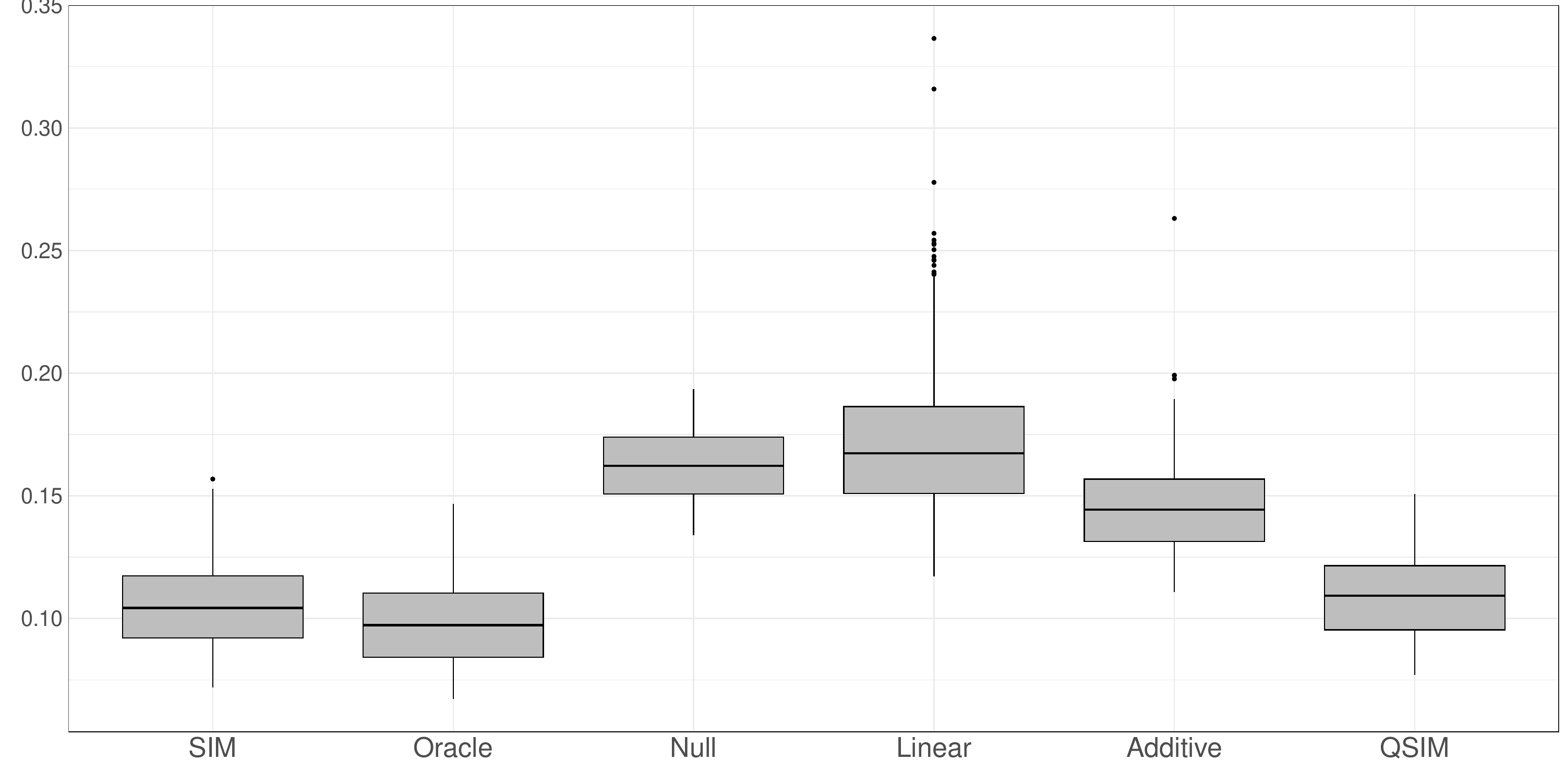}
\end{center}
\caption{Simulation results corresponding to Section 4.6.   
Boxplots of the ISE of the estimator for six estimators are presented.    \label{Result_Section46} }
\end{figure}

We calculated ISE as defined in an earlier section for each of the six competing estimators using 500 replications.
Figure \ref{Result_Section46} presents boxplots of the ISE for each estimator.
The results demonstrated that the distribution of the proposed estimator SIM closely aligns with that of Oracle, suggesting that the tuning parameter selection procedure described in Section 4.2 performed effectively.
By contrast, Null and Linear models are too simple to capture the nonlinear structure of $\alpha$ adequately.
Particularly, the Linear model exhibited large dispersion in its ISE, reflecting instability in its estimates.

The Additive model performance was superior to those of Null and Linear, but it was still inferior to that of the proposed estimator.
This finding is not surprising, given that the true model follows a single-index structure, which the Additive model cannot fully accommodate.

The QSIM also exhibited good performance.
However, in QSIM, the single-index parameter is estimated via linear quantile regression.
Because the true model is not linear in $\vec{X}$, this approach indicates the model misspecification, which appears to affect the overall estimation accuracy.
As a result, the ISE of QSIM was slightly larger than that of the proposed estimator SIM.

\section{Empirical Illustration}

The proposed single index model for EVI regression is applied to the motorcycle insurance claim data, which are available in the R package {\sf insuranceData} as {\sf dataOhlsson} (Ohlsson and Johansson, 2010).
This dataset comprises 64,548 motorcycle-related insurance records collected between 1994 and 1998 by the Swedish insurer Wasa.
Our primary objective is to model and predict the tail behavior of claim costs based on policyholder and policy characteristics.
However, approximately 99\% of the policies have zero claim cost. Only about 1\% led to positive claims.
To examine modeling large claim costs specifically, we restrict our analysis to the subset of 670 observations with positive claim costs.
Actually, if the data with zero claims were included, then the non-zero but minimum values of claim cost would be regarded as ^^ ^^ extreme values'' because they already lie in the top 1\% of the entire distribution. This would distort the interpretation of the tail behavior which we aim to model.
Daouia et al. (2022) and Zhang et al. (2024) also conducted statistical analyses of insurance data by removing observations with zero claim cost.
For this analysis, the response variable $Y$ is the claim cost. The covariates include seven standardized values of policy characteristics $\vec{X}=(X_1,\ldots,X_7)$, which are explained specifically in Table \ref{Ohlsson}.

\begin{table}
\caption{Descriptions of variables of motorcycle insurance claim data and corresponding single-index parameter estimates. \label{Ohlsson} }
\begin{tabular}{c|p{12cm}|>{\centering\arraybackslash}p{2.5cm}}
\hline
Symbol&Description& Single-index parameter\\
\hline
$Y$ & Claim cost\\
$X_1$ & Owner age, between 0 and 99 & 0.554\\
$X_2$ & Geographic zone numbered 1--7, in a standard classification of all Swedish parishes & 0.223\\
$X_3$ & MC class, a classification by the so-called EV ratio, defined as $({\rm Engine power in kW} \times
100) / ({\rm Vehicle weight in kilograms} + 75)$, rounded to the nearest lower integer. 75 kg denotes
the average driver weight. The EV ratios are divided into seven classes & -0.287\\
$X_4$ & Vehicle age, between 0 and 99 & 0.600\\
$X_5$ & Bonus class, taking values of 1--7. A new driver starts with bonus class 1. For each
claim-free year the bonus class is increased by 1. After the first claim the bonus is decreased
by 2. The driver can not return to class 7 with fewer than 6 consecutive claim free years & -0.435\\
$X_6$ & Number of policy years & -0.014\\
$X_7$ & Number of claims & 0.111\\
\hline
\end{tabular}
\end{table}

To apply our method to these data, we use $K=40$ equidistant knots on an interval 
\[
[-\min_i\|\vec{x}_i\|, \max_i\|\vec{x}_i\|].
\] 
Together with tuning parameters $(w_n,\lambda)$, we construct one estimator using the method presented in Sections 2.3 and 2.4. 
The top left panel shows the discrepancy measure $\{(w_s, D(w_s|\lambda_{s,cv}): s=1,\ldots,S\}$ with equidistant $S=300$ points for the 25\% quantile and the 90\% quantile of $Y$. 
The selected threshold value and the smoothing parameter were $w=51.93$. 
Then, the number of exceedances was $n_0=\sum_{i=1}^n I(Y_i>w)=115$. The exceedance rate was $115/670=0.172$. 

\begin{figure}
\begin{center}
\includegraphics[width=150mm,height=80mm]{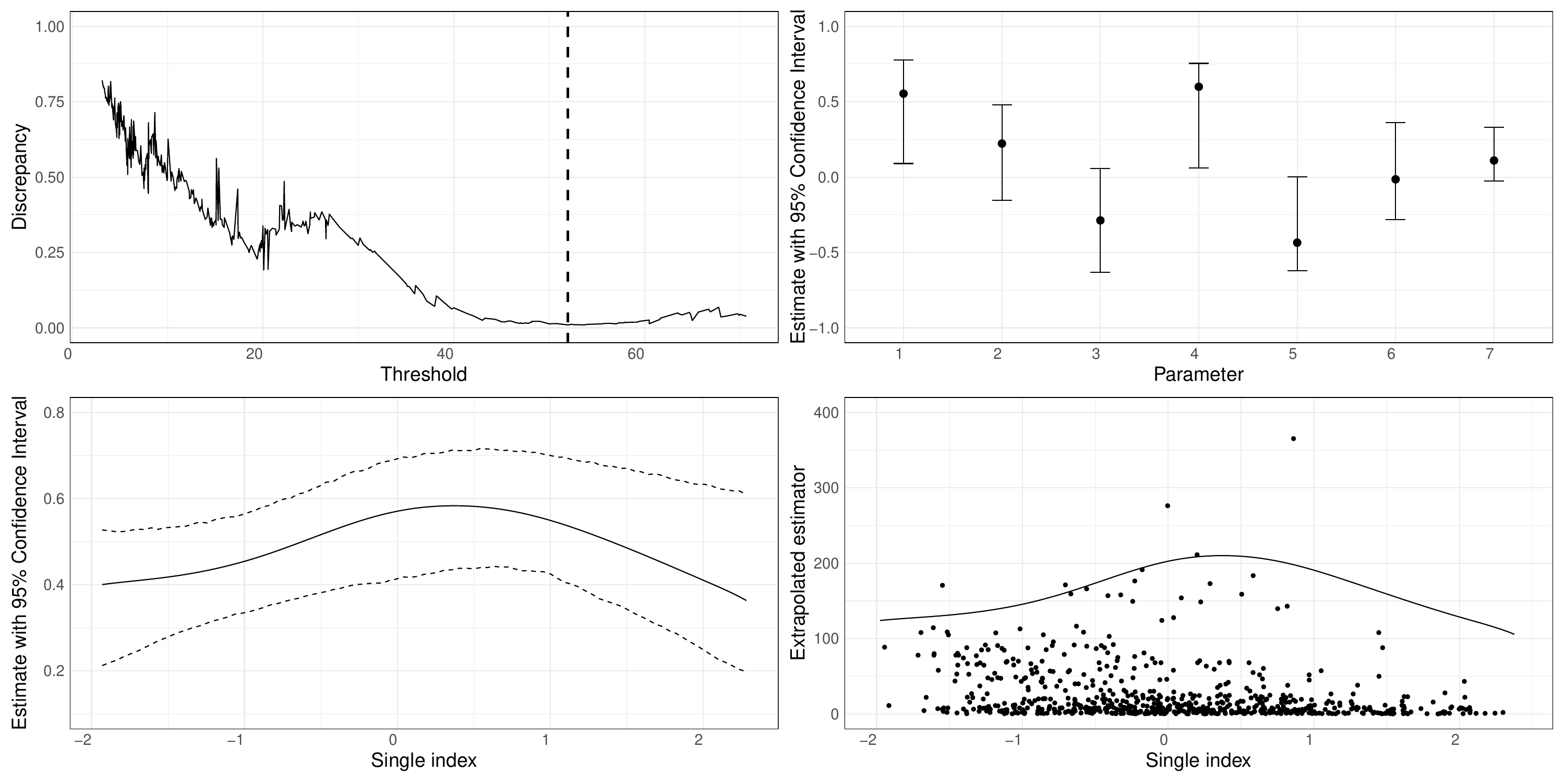}
\end{center}
\caption{Top left: Discrepancy measure over threshold (solid). The dashed line shows optimal $w_n$. 
Top right: Estimates and 95\% confidence interval of $\theta_1,\ldots,\theta_7$. 
Bottom left. Estimates and 95\% confidence interval of $\hat{\gamma}(z)$ for $z\in[-2,2.5]$. 
Bottom right: Extrapolated estimator of the 99\% conditional quantile of $Y$ given $\vec{X}^\top\hat{\vec{\theta}}$ with data $\{(\vec{X}_i^\top\hat{\vec{\theta}}, Y_i): i=1,\ldots,n, Y_i>0\}$.  \label{Result_Section5}}
\end{figure}

Using 115 exceedances, we estimate the single index parameter vector and nonlinear function. 
We also evaluate the estimation uncertainty from bootstrapping with 1000 replications. 
Then, for the fixed threshold $w=51.93$, we applied the bootstrap to the 115 data with $Y$ exceeding $w$ and calculate the estimator of $(\vec{\theta},\gamma)$.  
The estimator and 95\% confidence interval of each parameter are shown at the top right in Figure \ref{Result_Section5}. 
The point estimates of the single-index parameters are also listed in Table \ref{Ohlsson}.
It is apparent that the confidence interval was not symmetrical. 
However, that fact is not surprising because each $\theta_j$ is limited on $[-1,1]$ from the restriction $\|\hat{\vec{\theta}}\|=1$. 
The bottom left of Figure \ref{Result_Section5} shows $\hat{\gamma}(z)$ and the 95\% confidence interval at each $z\in[-2,2.5]$. 
From the result, at $z\in[0, 1]$, the large value of $\hat{\gamma}(z)$ was obtained. 

To obtain an easy interpretation of the behavior of the estimator of $\hat{\gamma}$, we constructed the estimator of the conditional quantile of $Y$ given $\vec{x}^\top\hat{\vec{\theta}}$ as  
\[
\tilde{Q}(\tau_E|\vec{x}^\top\hat{\vec{\theta}})=\left(\frac{n_0}{n(1-\tau_E)}\right)^{-\hat{\gamma}(\vec{x}^\top\hat{\vec{\theta}})}w_n,
\]
where $1-\tau_E < n_0/n$. 
The $\tilde{Q}$ is known as the extrapolated estimator of conditional quantile (Weismann 1978, Xu et al. 2020).
We set $\tau_E=0.99$. 
Because $n_0/n=0.17$, $1-\tau_E=0.01$ is much smaller than $n_0/n$. 
The behavior of the extrapolated estimator is described in the bottom right panel of Figure \ref{Result_Section5}. 
The result shows that the estimator $\tilde{Q}(\tau_E|\vec{x}^\top\hat{\vec{\theta}})$ exhibits a unimodal smooth curve that peaks near $\vec{x}^\top\hat{\vec{\theta}}=0.5$ and which decreases gradually as $\vec{x}^\top\hat{\vec{\theta}}$ moves away from it.
Furthermore, the resulting smooth curve implies that the influence of covariates on the extreme quantiles is stable and not erratic, which matches the expected pattern under our single index model with a smooth nonlinear function.

\section{Conclusion}

As described in this paper, we applied the single index model to the extreme value index (EVI) regression. 
Using the penalized maximum likelihood method for Pareto-type-tailed distribution approximation, we estimated the single index parameters and the one-dimensional nonlinear function included in the single index model. Additionally, we studied the asymptotic distribution and the rate of convergence of the proposed estimator. 
From these results, the single index model was confirmed as overcoming the curse of dimensionality. 
Simulation and empirical illustration help describe the efficiency of the proposed model. 

An important future task is variable selection when the dimension of covariates $p$ is quite large compared with sample size $n$ or sample size exceeding the threshold value. 
Nevertheless, no report of the relevant literature describes a result of sparse modelling or high-dimensional statistics in EVI regression. 
It would be interesting objective of future studies to investigate the hybrid method of high-dimensional statistics and extreme value theory. 

As described herein, we specifically examined only positive EVI and used the Pareto-type-tailed distribution. 
The single index model can be extended to general EVI including negative $\gamma$. 
For such cases, not only EVI but also the scale function (de Haan and Ferreira 2006) must be estimated. 
Although simultaneous estimation of the two target functions and establishing the asymptotic property of the estimator are quite difficult, exploring these aspects is extremely important.

\subsection*{Appendix A: $B$-spline basis} 

We now describe the definition and the property of the $B$-spline basis. 
Let $Z$ be a random variable with domain $[a,b]$.  
In this paper, we consider $Z=\vec{X}^\top\vec{\theta}$ for given $\vec{\theta}\in{\cal S}_+^{p-1}$. 
Again, we let $\vec{\kappa}=\{a=\kappa_0<\kappa_1<\ldots<\kappa_{K_0+1}=b\}, K_0>1$ be internal knots on an interval $[a, b]$.
Furthermore, let $\kappa_{-d+1}\leq\cdots\leq \kappa_{-1}<\kappa_0$ and $\kappa_{K_0+1}\leq\kappa_{K_0+2}\leq \cdots\leq \kappa_{K_0+d}$ be another set of knots.
For $j=-d+1,\ldots, K_0$ and $Z=z\in[a,b]$, let 
\[
\psi^{[0]}_j(z)=
\left\{
\begin{array}{cc}
1,\ & \kappa_j \leq z < \kappa_{j+1}\\
0, & otherwise
\end{array}
\right.
\]
be the firsr order ($d=1$) $B$-spline basis.
For $d>1$, the $d$th order $B$-spline basis can be defined recursively as 
\[
\psi_j^{[d]}(z)= \frac{z-\kappa_j}{\kappa_{j+d-1}-\kappa_j} \psi_j^{[d-1]}(z)+\frac{\kappa_{j+d}-z}{\kappa_{j+d}-\kappa_{j+1}} \psi_{j+1}^{[d-1]}(z),\ \ ^\forall z\in[a,b]. 
\]
For convenience, we treat $0/0=0$.
By the definition of $d$th order $B$-spline basis, we find that the $K=K_0+d-1$ basis function is used.
The scased $B$-spline basis are defined as 
\[
B_j^{[d]}(z)=\sqrt{K} \psi_j^{[d]}(z), j=-d+1,\ldots, K_0.
\]
The scaled $B$-spline basis is mathematically convenient than ordinary $B$-spline basis. 
Actually. since $\int \{\psi_j^{[d]}(z)\}^2 dz =O(K)$, we have $\int \{B_j^{[d]}(z)\}^2 dz =O(1)$ as $K\rightarrow\infty$.  
Although the normalized $B$-splines (Liu et al. 2011) are also useful, but in our model, centerization of $B$-splines is meaningless, and hence, only scale is adjusted. 
In particular, Lemma A.2 of Liu et al. (2011), which describes a key property of the scaled $B$-spline basis, is important and used in Appendices B--D below.

We see that the $m$th derivative of $B$-spline basis can be written by using $(d-m)$th degree $B$-spline basis. 
Actually, we see that for $m\geq 1$, 
\[
\frac{d^m \vec{B}^{[d]}(x)^\top\vec{b}}{d x^m }=\frac{d^m}{ dx^m}\sum_{j=1}^K B_j^{[d]}(x)b_j = \sum_{j=m+1}^K B_j^{[d-m]}(x) b_j^{(m)},
\]
where $\vec{b}=(b_1,\ldots,b_K)^\top\in\mathbb{R}^K$, 
\[
b_j^{(1)}= d \frac{b_j-b_{j-1}}{\kappa_{j+d}-\kappa_j}
\]
and 
\[
b_j^{(m)}= (d+1-m) \frac{b^{(m-1)}_j-b^{(m-1)}_{j-1}}{\kappa_{j+d+1-m}-\kappa_j}.
\]
This implies that the penalty term in (\ref{likelihood}) can be written as 
\[
\int \left[\frac{d^m \vec{B}^{[d]}(x)^\top\vec{b}}{d x^m }\right]^2 dx= \vec{b}^\top\Delta_{m,K}\vec{b},
\]
where $\Delta_{m,K}=D_{m,K}^\top R_m D_{m,K}$, $R_m$ is the $(K-m)$th square matrix having $(i,j)$-entry 
\[
\int B_i^{[d-m]}(x) B_j^{[d-m]}(x) dx
\]
and $D_{m,K}$ is the $(K-m)\times K$ matrix satisfying $\vec{b}^{(m)}=(b_{m+1}^{(m+1)},\ldots,b^{(m)}_{K})^\top=D_{m,K}\vec{b}$. 
If we use the equidistant knots $\kappa_j-\kappa_{j-1}=K^{-1}$, we obtain 
$D_m= K^{m} D^{diff}_{m,K}$, where $D^{diff}_{m,K}$ is the $m$th difference order matrix, which is defined as $D^{diff}_{m,K}=D^{diff}_{m-1,K-1}D^{diff}_{1,K}$ and $D^{diff}_{1,K}\vec{b}=(b_2-b_1,\ldots,b_{K}-b_{K-1})^\top$ (see Xiao 2019) . 
Consequently, the penalty term has the quadratic form with respect to $\vec{b}$.

\subsection*{Appendix B: Technical lemmas} 

We describe the technical lemmas used to the proof of theorems in Section 3.


\begin{lemma}\label{appspline}
Suppose that (C2) and (C5). Then, as $K\rightarrow\infty$,
\[
\sup_{z\in[a,b]} |\alpha_0(x)-\vec{B}(z)^\top\vec{b}_0|=O(K^{-d}).
\]
\end{lemma}

\begin{proof}[Proof of Lemma \ref{appspline}]

For simplicity, we write $X_0=\vec{X}^\top \vec{\theta}_0$ and $x_0=\vec{x}^\top \vec{\theta}_0$ for given $\vec{X}=\vec{x}$.
Define 
\begin{eqnarray*}
r_n(\vec{x})=\frac{(1/\gamma_0(\vec{x}^\top\vec{\theta})\beta(\vec{x})+1)^{-1}\ell_1(\vec{x}) w_n^{-\beta(\vec{x})}}{\ell_0(\vec{x})+\ell_1(\vec{x})w_n^{-\beta(\vec{x})}}.
\end{eqnarray*}
Then, from the definition of he Pareto-type tailed model (\ref{paretotail2}) with Hall class (\ref{slowly}), we have 
\begin{eqnarray}
&&E\left[\left.\frac{1}{\gamma_0(x_0)}\log\left(\frac{Y}{w_n}\right)\right|\vec{X}=\vec{x},Y>w_n\right]\nonumber\\
&&= \int_0^\infty P\left(\frac{1}{\gamma_0(x_0)}\log\left(\left.\frac{Y}{w_n}\right)>z\right|\vec{X}=\vec{x},Y>w_n\right)dz\nonumber\\
&&=
\int_0^\infty \frac{w_n^{-1/\gamma_0(x_0)}e^{-z}\{\ell_0(\vec{x})+\ell_1(\vec{x})w_n^{-\beta(\vec{x})}e^{-\gamma_0(x_0)\vec{\beta}(\vec{x})z}(1+o(1))\} }{w_n^{-1/\gamma_0(x_0)}\{\ell_0(\vec{x})+\ell_1(\vec{x})w_n^{-\beta(\vec{x})}(1+o(1))\}} dz\nonumber\\
&&=
\frac{\ell_0(\vec{x}) + (\gamma_0(x_0)\beta(\vec{x})+1)^{-1}\ell_1(\vec{x}) w_n^{-\beta(\vec{x})}(1+o(1))}{\ell_0(\vec{x}) + \ell_1(\vec{x}) w_n^{-\beta(\vec{x})}(1+o(1))}\nonumber\\
&&=
1+ r_n(\vec{x})(1+o(1)). \label{BiasSecond}
\end{eqnarray}
This implies that 
\begin{eqnarray*}
E\left[\left. \log\left(\frac{Y}{w_n}\right)\right| Y>w_n, \vec{X}=\vec{x}\right]
= \gamma_0(x_0)\{1+r_n(\vec{x})(1+o(1))\}.
\end{eqnarray*}
Since $L(\vec{b})$ is convex function, the minimizer of $\vec{b}_0$ is unique and this is the solution of 
\begin{eqnarray*}
\frac{\partial}{\partial \vec{b}} L(\vec{b})
&=&
E\left[\left.\left\{\exp[\vec{B}(X_0)^\top\vec{b}_0]\log\left(\frac{Y}{w_n}\right)-1\right\}\vec{B}(X_0)\right| Y>w_n\right]\\
&=&
E\left[\left\{\exp[\vec{B}(X_0)^\top\vec{b}_0-\alpha_0(X_0)]\{1+r_n(\vec{X})(1+o(1))\}-1\right\}\vec{B}(X_0)\right]\\
&=&
\vec{0}.
\end{eqnarray*}
Thus, if $\partial L(\vec{b}_0)/\partial \vec{b}=\vec{0}$, $\vec{B}(X_0)^\top\vec{b}_0-\alpha_0(X_0)$ must be as small as possible. 
Meanwhile from Barrow and Smith (1987), for $\alpha_0\in{\cal C}^q[a,b]$ with $d\leq q$, there exists $\vec{b}^*\in\mathbb{R}^K$ such that $\sup_{z\in[a,b]}|\alpha_0(z)-\vec{B}(z)^\top\vec{b}^*|=O(K^{-d})$. 
Therefore, if $K^d |\alpha_0(x)-\vec{B}(z)^\top\vec{b}_0|\rightarrow \infty$, we obtain $L(\vec{b}^*) < L(\vec{b}_0)$, which contradict the fact that $\vec{b}_0$ is the minimizer of $L(\vec{b})$. 
This implies that $ |\alpha_0(x)-\vec{B}(z)^\top\vec{b}_0|=O(K^{-d})$ for any $z\in[a,b]$.
Thus, Lemma \ref{appspline} was proven.
\end{proof}

Here, for a square matrix $A$, let $\rho_{min}(A)$ and $\rho_{max}(A)$ be the minimum and maximum eigen value of $A$, respectively. 
We define
\begin{eqnarray*}
\Sigma
=
\left[
\begin{array}{cc}
\Sigma_{b,b} &\Sigma_{b,\phi}\\
\Sigma_{\phi,b} &\Sigma_{\phi,\phi}\\
\end{array}
\right],
\end{eqnarray*}
where
\begin{eqnarray*}
\Sigma_{b,b}
&=& E[P(Y>w_n|\vec{X})\vec{B}(\vec{X}^\top\vec{\theta}_0)\vec{B}(\vec{X}^\top\vec{\theta}_0)^\top] +\lambda\Delta_{m,k},\\
\Sigma_{b,\phi}
&=&
E[P(Y>w_n|\vec{X}\alpha^{(1)}(\vec{X}^\top\vec{\theta}_0) \vec{B}(\vec{X}^\top\vec{\theta}_0)\vec{X}^\top J_1(\vec{\phi})^\top],
\end{eqnarray*}
$\Sigma_{\phi,b}=\Sigma_{b,\phi}^\top$,
\[
\Sigma_{\phi,\phi}=
E\left[P(Y>w_n|\vec{X})\{\alpha^{(1)}(\vec{X}^\top\vec{\theta}_0)\}^2J_1(\vec{\phi}_0) XX^\top J_1(\vec{\phi}_0)^\top\right],
\]
$p$-idendity matrix $I_p$,  $(p-1)\times p$ matrix $J_1(\vec{\phi})=[\vec{\phi}/\sqrt{1-\|\vec{\phi}\|}\ I_{p-1}]$ and $\alpha_0^{(j)}(z)= d^j \alpha_0(z)/d z^j$.
From Lemma \ref{hessian}, we have 
\begin{eqnarray*}
E\left[
\begin{array}{cc}
\frac{\partial^2 \ell_n(\vec{b}_0,\vec{\phi}_0)}
{\partial \vec{b}\partial\vec{b}^\top}
&
\frac{\partial^2 \ell_n(\vec{b}_0,\vec{\phi}_0)}
{\partial \vec{b}\partial\vec{\phi}^\top}
\\
\frac{\partial^2 \ell_n(\vec{b}_0,\vec{\phi}_0)}
{\partial\vec{\phi} \partial \vec{b}^\top} & 
\frac{\partial \ell_n(\vec{b}_0,\vec{\phi}_0)}
{\partial \vec{\phi} \partial \vec{\phi}}
\end{array}
\right]
=\Sigma(1+o(1)).
\end{eqnarray*}
That is, the matrix $\Sigma$ is the Hessian matrix of objective penalized log-likelihood function.

\begin{lemma}\label{EigenHessian}
Suppose that (C1)--(C6). 
Then, there exist constants $C_*>0$ and $C^*>0$ such that $\rho_{min}(\Sigma) \geq C_* \tau_n$ and $\rho_{max}(\Sigma)\leq C^*\tau_n$ as $n\rightarrow\infty$.
\end{lemma}

\begin{proof}[Proof of Lemma \ref{EigenHessian}]

Under (C1), $\Sigma$ is non-singular, and hence it is sufficient to show $\rho_{min}(\Sigma)=O(\tau_n)$ and $\rho_{max}(\Sigma)=O(\tau_n)$. 
Let $\vec{u}\in\mathbb{R}^{K+p-1}-\{\vec{0}\}$ with $\|\vec{u}\|=1$. 
We write $\vec{u}=(\vec{u}_b^\top, \vec{u}_\phi^\top)^\top$, where $\vec{u}_b\in\mathbb{R}^K$ and $\vec{u}_\phi\in\mathbb{R}^{p-1}$. 
We note that for $\vec{u}_b=(u_{1,b},\dots, u_{K,b})^\top$, $\max_j|u_{j,b}|=O(K^{-1/2})$ since $\|\vec{u}_b\|^2 <1$. 
For 
\[
\vec{u}^\top \Sigma \vec{u}=\vec{u}_b^\top \Sigma_{b,b}\vec{u}_b+2\vec{u}_b^\top \Sigma_{b,\phi}\vec{u}_\phi + \vec{u}_\phi^\top \Sigma_{\phi,\phi}\vec{u}_\phi,
\]
we first consder $\vec{u}_\phi^\top \Sigma_{\phi,\phi}\vec{u}_\phi$. 
From the definition of $\Sigma_{\phi,\phi}$ and the mean value theorem for integrals, there exists $\vec{x}_*$ such that 
\begin{eqnarray*}
\vec{u}_\phi^\top \Sigma_{\phi,\phi}\vec{u}_\phi&=&
E\left[P(Y>w_n|\vec{X})\{\alpha^{(1)}(\vec{X}^\top\vec{\theta}_0)\}^2\vec{u}_\phi^\top J_1(\vec{\phi}_0) XX^\top J_1(\vec{\phi}_0)^\top\right]\\
&=& E\left[P(Y>w_n|\vec{X})\right] \{\alpha^{(1)}(\vec{x}_*^\top\vec{\theta}_0)\}^2J_1(\vec{\phi}_0) \vec{x}_*\vec{x}_*^\top J_1(\vec{\phi}_0)^\top \vec{u}_\phi\\
&=&\tau_n \{\alpha^{(1)}(\vec{x}_*^\top\vec{\theta}_0)\}^2 \{\vec{u}_\phi^\top J_1(\vec{\phi}_0) \vec{x}_*\}^2\\
&=&O(\tau_n).
\end{eqnarray*}
Next, we evaluate 
\[
\vec{u}_b^\top\Sigma_{b,\phi}\vec{u}_\phi
=
E[P(Y>w_n|\vec{X})\alpha^{(1)}(\vec{X}^\top\vec{\theta}_0) \vec{u}_b^\top\vec{B}(\vec{X}^\top\vec{\theta}_0)\vec{X}^\top J_1(\vec{\phi})^\top\vec{u}_\phi].
\]
Under (C1), for any $\vec{x}\in{\cal X}$, $\vec{X}^\top J_1(\vec{\phi})^\top\vec{u}_\phi=O(1)$. 
Next, from Appendix A, for any $z= \vec{x}^\top\vec{\theta}_0$ with $\vec{x}\in{\cal X}$, there exists $j^*$ such that 
\begin{eqnarray}
\vec{u}_b^\top\vec{B}(z)=\sum_{j=1}^K B_j(z)u_{j,b} = \sum_{j=j^*}^{j^* +d} B_j(z)u_{j,b}  \label{BspPro}
\end{eqnarray}
and $B_j(z)=0$ for $j<j^*$ and $j>j^*+d$. 
Since $B_j(z)=O(\sqrt{K})$ and $u_{j,b}=O(K^{-1/2})$, we obtain $|\vec{u}_b^\top\vec{B}(z)|=O(1)$. 
Therefore, from mean value theorem for integrals, there exists $\vec{x}_*\in{\cal X}$ such that 
\begin{eqnarray*}
\vec{u}_b^\top\Sigma_{b,\phi}\vec{u}_\phi
&=&
E[P(Y>w_n|\vec{X})]\alpha^{(1)}(\vec{x}_*^\top\vec{\theta}_0) \vec{u}_b^\top\vec{B}(\vec{x}_*^\top\vec{\theta}_0)\vec{x}_*^\top J_1(\vec{\phi})^\top\vec{u}_\phi =O(\tau_n). 
\end{eqnarray*}
Lastly, we consider 
\begin{eqnarray*}
\vec{u}_b^\top\Sigma_{b,b}\vec{u}_b
=E[P(Y>w_n|\vec{X})\vec{u}_b^\top\vec{B}(\vec{X}^\top\vec{\theta}_0)\vec{B}(\vec{X}^\top\vec{\theta}_0)^\top\vec{u}_b] +\lambda \vec{u}_b^\top\Delta_{m,k}\vec{u}_b.
\end{eqnarray*}
Similar to (\ref{BspPro}), for any $z\in[a,b]$, there exists $j^*$ such that
\begin{eqnarray*}
\{\vec{u}_b^\top\vec{B}(z)\}^2=\left\{\sum_{j=1}^K B_j(z)u_{j,b} \right\}^2= \left\{\sum_{j=j^*}^{j^* +d} B_j(z)u_{j,b}\right\}^2 
\end{eqnarray*}
and $B_j(z)=0$ for $j<j^*$ and $j>j^*+d$. 
From $B_j(z)=O(K^{1/2})$ and $u_{j,b}=O(K^{-1/2})$, we obtain $\{\vec{u}_b^\top\vec{B}(z)\}^2=O(1)$, which is standard property of scaled $B$-spline model. 
Therefore, mean value of theorem for integrals yelds that there exists $\vec{x}_*\in{\cal X}$ such that 
\[
E[P(Y>w_n|\vec{X})\vec{u}_b^\top\vec{B}(\vec{X}^\top\vec{\theta}_0)\vec{B}(\vec{X}^\top\vec{\theta}_0)^\top\vec{u}_b] =E[P(Y>w_n|\vec{X})]\{\vec{u}_b^\top\vec{B}(\vec{x}_*^\top\vec{\theta}_0)\}^2=O(\tau_n).
\]
Next, from the Proposition 4.2 of Xiao (2019) and (C6), we have 
\[
0\leq \lambda \vec{u}_b^\top\Delta_{m,K}\vec{u}_b=O(\lambda K^{2m}) =O(\tau_n).
\]
Conseqently, Lemma \ref{EigenHessian} holds. 

\end{proof}

We next show the expectation of gradient of $U_n$. 
\begin{lemma}\label{gradient}
Suppose that (C1)--(C6). As $n\rightarrow\infty$, 
\[
\left\|E\left[
\frac{\partial \ell_n(\vec{b}_0,\vec{\phi}_0)}
{\partial \vec{b}}\right]\right\|^2
\leq O(\tau_n^{2+2\beta_{inf}}K) + O(\tau_n\lambda K)
\]
and 
\[
\left\|E\left[\frac{\partial \ell_n(\vec{b}_0,\vec{\phi}_0)}
{\partial \vec{\phi}}
\right]\right\|^2
\leq O(\tau_n^{2+2\beta_{inf}})
\]
\end{lemma}

\begin{proof}[Proof of Lemma \ref{gradient}]

For simplicity, we write $X_{0i}=\vec{X}_i^\top\vec{\theta}_0=\vec{X}_i^\top\vec{\theta}(\vec{\phi}_0)$ and $X_0=\vec{X}^\top\vec{\theta}(\vec{\phi}_0)$. 
We then obtain
\begin{eqnarray*}
\frac{\partial \ell_n(\vec{b}_0,\vec{\phi}_0)}
{\partial \vec{b}}
= \frac{1}{n}\sum_{i=1}^n \left\{\exp[\vec{B}(X_{0i})^\top\vec{b}_0]\log\left(\frac{Y_i}{w_n}\right)-1 \right\}\vec{B}(X_{0i})I(Y_i>w_n) +\lambda\Delta_{m,K}\vec{b}_0.
\end{eqnarray*}
Under (C3) and (C5), we obtain $K^{-q}/\inf_{\vec{x}\in{\cal X}}r_n(\vec{x})\rightarrow 0$. 
From this, Lemma 1 and (\ref{BiasSecond}), we have 
\begin{eqnarray*}
E\left[\frac{\partial \ell_n(\vec{b}_0,\vec{\phi}_0)}
{\partial \vec{b}}\right]
&=& E\left[P(Y>w_n|\vec{X})\left.\left\{\exp[\vec{B}(X_{0i})^\top\vec{b}_0]\log\left(\frac{Y_i}{w_n}\right)-1 \right\}\vec{B}(X_{0i})\right| Y>w_n\right] +\lambda\Delta_{m,K}\vec{b}_0\\
&=&
E[P(Y>w_n|\vec{X})r_n(\vec{X})\vec{B}(X_0)](1+o(1))+ \lambda\Delta_{m,K}\vec{b}_0.
\end{eqnarray*}
Since $E[\partial \ell_n(\vec{b}_0,\vec{\phi}_0)/\partial \vec{b}]$ is $K$-vector, the asymptotic order of the squared norm of this is similar to that of
\[
 O(K) \times \rho_{\max}\left(E\left[\frac{\partial \ell_n(\vec{b}_0,\vec{\phi}_0)}
{\partial \vec{b}}\right]E\left[\frac{\partial \ell_n(\vec{b}_0,\vec{\phi}_0)}
{\partial \vec{b}}\right]^\top\right).
\]
We aim is to show 
\[
\rho_{\max}\left(E\left[\frac{\partial \ell_n(\vec{b}_0,\vec{\phi}_0)}
{\partial \vec{b}}\right]E\left[\frac{\partial \ell_n(\vec{b}_0,\vec{\phi}_0)}
{\partial \vec{b}}\right]^\top\right)\leq  O(\tau_n^{2+2\beta_{inf}})+O(\tau_n\lambda).
\]
To this ends, we consider that for $\vec{u}\in\mathbb{R}^K$ with $\|\vec{u}\|=1$,
\begin{eqnarray*}
\vec{u}^\top E\left[\frac{\partial \ell_n(\vec{b}_0,\vec{\phi}_0)}
{\partial \vec{b}}\right]
=
E[P(Y>w_n|\vec{X})r_n(\vec{X}) \vec{u}^\top\vec{B}(X_0)] + \lambda \vec{u}^\top \Delta_{m,K}\vec{b}_0
\end{eqnarray*}
The mean value theorem for integrals yeilds that there exists $z^*\in[a,b]$ such that 
\[
E[P(Y>w_n|\vec{X})r_n(\vec{X}) \vec{u}^\top\vec{B}(X_0)] =\vec{u}^\top\vec{B}(z^*) E[P(Y>w_n|\vec{X})r_n(\vec{X})].
\]
Similar to the proof of Lemma \ref{EigenHessian}, we have $|\vec{u}^\top\vec{B}(z^*)|=O(1)$. 
Meanwhile, from the definition (\ref{paretotail}) and (\ref{slowly}), we have $P(Y>w_n|\vec{x})\approx \ell_0(\vec{x})w_n^{-1/\gamma(\vec{x}^\top\vec{\theta}_0)}$. 
Therefore, we obtain 
\[
|r_n(\vec{x})| \leq C^*  \ell_0(\vec{x})^{1/\gamma(\vec{x}^\top\vec{\theta}_0)\beta(\vec{x})}w_n^{-\beta(\vec{x})} \leq C^* P(Y>w_n|\vec{x})^{\gamma(\vec{x}^\top\vec{\theta}_0)\beta(\vec{x})} \leq C^* P(Y>w_n|\vec{x})^{\beta_{inf}},
\]
where $C^*$ is a constant satisfying 
\[
\left|\frac{(1/\gamma_0(\vec{x}^\top\vec{\theta})\beta(\vec{x})+1)^{-1}\ell_1(\vec{x})}{\ell_0(\vec{x})^{1/\gamma(\vec{x}^\top\vec{\theta}_0)\beta(\vec{x}) }(\ell_0(\vec{x})+\ell_1(\vec{x})w_n^{-\beta(\vec{x})})}\right| \leq C^*.
\]
Note that the finiteness of $C^*$ can be guaranteed by (C3). 
Next, from (C5) and the definition of $\Delta_{m,K}$ in Appendix A, 
\[
\lambda\vec{u}^\top\Delta_{m,K} \vec{b}_0 =\lambda\vec{u}^\top D_{m,K}^\top \int \vec{B}^{[d-m]}(x) \alpha_0^{(m)}(x) dx(1+o(1))
\]
and each element of $\Delta_{m,K}$ is $O(K^m)$.
Since $\|\vec{u}\|=1$, each element of $\vec{u}$ has $O(K^{-1/2})$. Appendix A
In addition, the property of scaled $B$-spline model shows $\int \vec{B}^{[d-m]}(x) \alpha_0^{(m)}(x) dx=O(K^{-1/2})$. 
Together with the fact that $\Delta_{m,K}$ is band matrix, we obtain $\lambda\vec{u}^\top\Delta_{m,K} \vec{b}_0 =O(\lambda K^m)$.
Thus, we obtain 
\[
\{E[P(Y>w_n|\vec{X})r_n(\vec{X}) \vec{u}^\top\vec{B}(X_0)] \}^2 \leq O(\tau_n^{2+2\beta_{inf}})+O(\lambda^2 K^{2m}). 
\]
Under (C6), we have $\lambda K^{2m}/\tau_n=O(1)$, which implies that $O(\lambda^2K^{2m})=O(\tau_n\lambda)$

Next, we consider $\partial U_n(\vec{b}_0,\vec{\phi}_0)/\partial \vec{\phi}$. 
From the definition of $\vec{\theta}(\vec{\phi}_0)$ and Lemma 1, we have 
\[
\frac{\partial}{\partial \vec{\phi}} \vec{B}(\vec{X}^\top \vec{\theta}(\vec{\phi}))^\top\vec{b}_0
= J_1(\vec{\phi})\vec{X} \alpha^{(1)}(\vec{X}^\top \vec{\theta}(\vec{\phi}))(1+o(1)).
\]
This and (\ref{BiasSecond}) imply 
\begin{eqnarray*}
E\left[\frac{\partial \ell_n(\vec{b}_0,\vec{\phi}_0)}
{\partial \vec{\phi}}\right]
&=&
E\left[P(Y>w_n|\vec{X}) \alpha_0^{(1)}(X_0)J(\vec{\phi}_0)\vec{X}\left.\left\{\exp[\vec{B}(X_0)^\top\vec{b}_0]\log\left(\frac{Y_i}{w_n}\right)-1\right\}\right| Y>w_n \right]\\
&=&
E\left[P(Y>w_n|\vec{X}) \alpha_0^{(1)}(X_0)J(\vec{\phi}_0)\vec{X}r_n(\vec{X})(1+o(1))\right].
\end{eqnarray*}
Thus, Lemma \ref{gradient} was shown.
\end{proof}

\begin{lemma}\label{hessian}
Suppose that (C1)--(C6). As $n\rightarrow\infty$, 
\begin{eqnarray*}
E\left[
\begin{array}{cc}
\frac{\partial^2 \ell_n(\vec{b}_0,\vec{\phi}_0)}
{\partial \vec{b}\partial\vec{b}^\top}
&
\frac{\partial^2 \ell_n(\vec{b}_0,\vec{\phi}_0)}
{\partial \vec{b}\partial\vec{\phi}^\top}
\\
\frac{\partial^2 \ell_n(\vec{b}_0,\vec{\phi}_0)}
{\partial\vec{\phi} \partial \vec{b}^\top} & 
\frac{\partial \ell_n(\vec{b}_0,\vec{\phi}_0)}
{\partial \vec{\phi} \partial \vec{\phi}}
\end{array}
\right]
=\Sigma(1+o(1)).
\end{eqnarray*}
\end{lemma}

\begin{proof}[Proof of Lemma \ref{hessian}]
Similar to proof of Lemma \ref{gradient}, we write $X_{0i}=\vec{X}_i^\top\vec{\theta}_0$ and $X_0=\vec{X}^\top\vec{\theta}_0$. 
We note that $\partial U_n(\vec{b},\vec{\phi})/\partial \vec{b}$ and $\partial U_n(\vec{b},\vec{\phi})/\partial \vec{\phi}$ are already shown in the proof of Lemma \ref{gradient}.
We first obtain
\begin{eqnarray*}
E\left[\frac{\partial^2 \ell_n(\vec{b}_0,\vec{\phi}_0)}
{\partial \vec{b}\partial\vec{b}^\top}\right]
&=&
E\left[
P(Y>w_n|\vec{X})\vec{B}(X_0)\vec{B}(\vec{X}^\top\vec{\theta}(\vec{\phi}))^\top \exp[\vec{B}(X_0)^\top\vec{b}_0]\log\left(\frac{Y_i}{w_n}\right)
|Y>w_n\right]\\
&=&\Sigma_{b,b}(1+o(1)).
\end{eqnarray*}
Next, we have 
\begin{eqnarray*}
&&E\left[\frac{\partial^2 \ell_n(\vec{b}_0,\vec{\phi}_0)}
{\partial \vec{b}\partial\vec{\phi}^\top}\right]\\
&&=
E\left[\left.\frac{\partial}{\partial \vec{\phi}^\top}P(Y>w_n|\vec{X})\left\{\exp[\vec{B}(X_0)^\top\vec{b}_0]\log\left(\frac{Y_i}{w_n}\right)-1 \right\} \vec{B}(\vec{X}^\top\vec{\theta}(\vec{\phi})) \right|_{\vec{\phi}=\vec{\phi}_0} | Y>w_n\right]\\
&&=
E\left[P(Y>w_n|\vec{X})\alpha^{(1)}_0(X_0)\vec{B}(\vec{X}^\top\vec{\theta}(\vec{\phi}))\vec{X}^\top J_1(\vec{\phi}_0)^\top\right](1+o(1))\\
&&\quad + 
E\left[P(Y>w_n|\vec{X})\left\{\exp[\vec{B}(X_0)^\top\vec{b}_0]\log\left(\frac{Y_i}{w_n}\right)-1 \right\} \frac{\partial B(\vec{X}^\top\vec{\theta}(\vec{\phi}_0))}{\partial \vec{\phi}}| Y>w_n \right]\\
&&=
\Sigma_{b,\phi}(1+o(1))
+
E\left[P(Y>w_n|\vec{X})r_n(\vec{X})\frac{\partial B(\vec{X}^\top\vec{\theta}(\vec{\phi}_0))}{\partial \vec{\phi}}| Y>w_n \right](1+o(1)).
\end{eqnarray*}
Since the asymptotic order of $\vec{B}(z)$ and $\partial \vec{B}(z)/\partial z$ are similar (see, de Boor 2001), we obtain 
\[
\left|E\left[P(Y>w_n|\vec{X})r_n(\vec{X})\frac{\partial B(\vec{X}^\top\vec{\theta}(\vec{\phi}_0))}{\partial \vec{\phi}}| Y>w_n \right]\right|\leq
O(\tau_n^{1+\beta_{inf}}),
\] 
which is negrigible order compared with $\Sigma_{b,\phi}=O(\tau_n)$. 

Let $(p-1)$-matrix $J_2(\vec{\phi})=(J_{2,i,j})$, where $\vec{\phi}=(\phi_1,\ldots,\phi_{p-1})$ and
\[
J_{2,i,j}=\frac{1}{\sqrt{1-\|\vec{\phi}\|^2}} \left(1- \frac{1}{2}\frac{\phi_i\phi_j}{1-\|\vec{\phi}\|^2}\right), \ \ i,j=1,\ldots,p-1.
\]
Then, we have 
\[
\frac{\partial}{\partial \vec{\phi}^\top} J_1(\vec{\phi}) = J_2(\vec{\phi}).
\]
We then obtain
\begin{eqnarray*}
&&E\left[\frac{\partial^2 \ell_n(\vec{b}_0,\vec{\phi}_0)}
{\partial \vec{\phi}\partial\vec{\phi}^\top}\right]\\
&&=
E\left[\left.\frac{\partial }{\partial \vec{\phi}} P(Y>w_n|\vec{X})\alpha_0^{(1)}(\vec{X}^\top\vec{\theta}(\vec{\phi}))J_1(\vec{\phi})\vec{X}\left\{\exp[\alpha_0(\vec{X}^\top\vec{\theta}(\vec{\phi}))]\log\left(\frac{Y_i}{w_n}\right) -1\right\}\right|_{\vec{\phi}=\vec{\phi}_0} | Y>w_n\right]\\
&&\times (1+o(1))\\
&&=
E\left[P(Y>w_n|\vec{X})\alpha_0^{(2)}(\vec{X}^\top\vec{\theta}(\vec{\phi}_0))J_1(\vec{\phi}_0)\vec{X}\vec{X}^\top J_1(\vec{\phi}_0)r_n(\vec{X})\right](1+o(1))\\
&&\quad + 
E\left[P(Y>w_n|\vec{X})\alpha_0^{(1)}(\vec{X}^\top\vec{\theta}(\vec{\phi}_0))J_2(\vec{\phi}_0)X_1 r_n(\vec{X})\right](1+o(1))\\
&&\quad\quad + E\left[P(Y>w_n|\vec{X})\{\alpha_0^{(1)}(X_0)\}^2J_1(\vec{\phi}_0)\vec{X}^\top\vec{X}J_1(\vec{\phi})^\top\right](1+o(1)).
\end{eqnarray*}
Under (C5), we have 
\[
\left|E\left[P(Y>w_n|\vec{X})\alpha_0^{(2)}(\vec{X}^\top\vec{\theta}(\vec{\phi}_0))J_1(\vec{\phi})\vec{X}\vec{X}^\top J_1(\vec{\phi}_0)r_n(\vec{X})\right]\right|\leq O(\tau_n^{1+\beta_{inf}})
\]
and 
\[
\left|E\left[P(Y>w_n|\vec{X})\alpha_0^{(1)}(\vec{X}^\top\vec{\theta}(\vec{\phi}_0))J_2(\vec{\phi}_0)X_1 r_n(\vec{X})\right]\right|\leq O(\tau_n^{1+\beta_{inf}}).
\]
These orders are smaller order than 
\[
\Sigma_{\phi,\phi} =E\left[P(Y>w_n|\vec{X})\{\alpha_0^{(1)}(X_0)\}^2J_1(\vec{\phi}_0)\vec{X}^\top\vec{X}J_1(\vec{\phi}_0)^\top\right]=O(\tau_n).
\]
Thus, Lemma \ref{hessian} was proven. 
\end{proof}

\begin{lemma}\label{Consistency}
Suppose that (C1)--(C6). Then, as $n\rightarrow\infty$,
\[
\|\hat{\vec{b}}-\vec{b}_0\| + \|\hat{\vec{\phi}}-\vec{\phi}_0\| \stackrel{P}{\to} 0.
\]
\end{lemma}

\begin{proof}[Proof of Lemma \ref{Consistency}]

Let 
\[
L_0(\vec{b}, \vec{\phi})=E\left[\left\{\exp[\vec{B}(\vec{X}^\top\vec{\theta}(\vec{\phi}))^\top\vec{b}]\log\left(\frac{Y}{w_n}\right)-\vec{B}(\vec{X}^\top\vec{\theta}(\vec{\phi}))^\top\vec{b}\right\}I(Y>w_n)\right]
\]
and
\begin{eqnarray*}
L(\vec{b},\vec{\phi})&=&\ell_n(\vec{b},\vec{\theta}(\vec{\phi})|\lambda)\\
&=&\frac{1}{n}\sum_{i=1}^n \left\{\exp[\vec{B}(\vec{X}_i^\top\vec{\theta}(\vec{\phi}))^\top\vec{b}]\log\left(\frac{Y_i}{w_n}\right)-\vec{B}(\vec{X}_i^\top\vec{\theta}(\vec{\phi}))^\top\vec{b}\right\}I(Y_i>w_n)\\
&&+\frac{\lambda}{2}\int_a^b \left\{\frac{d^m}{dx^m}\vec{B}(z)^\top\vec{b} \right\}^2dz.
\end{eqnarray*}
We note that $L$ and $L_0$ are strictly convex functions. 
Therefore, $(\vec{b}_0,\vec{\phi}_0)=\argmin_{\vec{b},\vec{\phi}} L_0(\vec{b},\vec{\phi})$ and $(\hat{\vec{b}},\hat{\vec{\phi}})=\argmin_{\vec{b},\vec{\phi}} L(\vec{b}, \vec{\phi})$ are uniquely defined. 
Define $\eta(\vec{b},\vec{\phi})=\|\vec{b}-\vec{b}_0\|^2+\|\vec{\phi}-\vec{\phi}_0\|^2$. 
From Lemma 2 of Hijort and Pollard (1993), for any $\varepsilon>0$,
\begin{eqnarray*}
&&P(\|\hat{\vec{b}}-\vec{b}_0\|^2+\|\hat{\vec{\phi}}-\vec{\phi}_0\|^2>\varepsilon^2)\\
&&\leq P\left(\sup_{\eta(\vec{b},\vec{\phi})\leq \varepsilon^2}|L(\vec{b},\vec{\phi})-L_0(\vec{b},\vec{\phi})|\geq 2^{-1}\inf_{\eta(\vec{b},\vec{\phi})= \varepsilon^2}|L_0(\vec{b},\vec{\phi})-L_0(\vec{b}_0,\vec{\phi}_0)| \right).
\end{eqnarray*}
We now consider the vector $\vec{b}\in\mathbb{R}^K$ and $\vec{\phi}\in\mathbb{R}^p$ satisfying $\eta(\vec{b},\vec{\phi})= \varepsilon^2$.
We then write $\vec{b}=\vec{b}_0+ \varepsilon^2 \vec{u}_b, \vec{u}_b\in\mathbb{R}^{K}$ and $\vec{\phi}=\vec{\phi}_0+ \varepsilon^2 \vec{u}_\phi, \vec{u}_\phi\in\mathbb{R}^{p-1}$, where $\|\vec{u}_b\|^2+\|\vec{u}_\phi\|^2=1$. 
Since the hessian of $L_0$ is continuous with respect to $(\vec{b},\vec{\phi})$, by the Taylor's theorem, we obtain
\begin{eqnarray*}
&&L_0(\vec{b},\vec{\phi})-L_0(\vec{b}_0,\vec{\phi}_0)\\
&&=
\varepsilon^2\left\{\frac{\partial L_0(\vec{b}_0,\vec{\phi}_0)}{\partial \vec{b}^\top}\vec{u}_b + \frac{\partial L_0(\vec{b}_0,\vec{\phi}_0)}{\partial \vec{\phi}^\top}\vec{u}_\phi\right\}\\
&&+
\varepsilon^4\left\{\vec{u}_b^\top\frac{\partial^2 L_0(\vec{b}_0,\vec{\phi}_0)}{\partial \vec{b}\partial \vec{b}^\top}\vec{u}_b + 
2\vec{u}_b^\top\frac{\partial^2 L_0(\vec{b}_0,\vec{\phi}_0)}{\partial \vec{b}\partial \vec{\phi}^\top}\vec{u}_\phi
+
\vec{u}_\phi^\top\frac{\partial^2 L_0(\vec{b}_0,\vec{\phi}_0)}{\partial \vec{\phi}\partial \vec{\phi}^\top}\vec{u}_\phi\right\}(1+o(1)).
\end{eqnarray*}
By the definition of $L_0$, $\partial L_0(\vec{b}_0,\vec{\phi}_0)/\partial \vec{b}=\vec{0}$ and $\partial L_0(\vec{b}_0,\vec{\phi}_0)/\partial \vec{\phi}=\vec{0}$. 
Since $\|\vec{u}_b\|<1$ and $\|\vec{u}_\phi\|<1$, from Lemma \ref{EigenHessian}, there exists a constant $c^*>0$ such that 
\[
\left|\vec{u}_b^\top\frac{\partial^2 L_0(\vec{b}_0,\vec{\phi}_0)}{\partial \vec{b}\partial \vec{b}^\top}\vec{u}_b + 
2\vec{u}_b^\top\frac{\partial^2 L_0(\vec{b}_0,\vec{\phi}_0)}{\partial \vec{b}\partial \vec{\phi}^\top}\vec{u}_\phi
+
\vec{u}_\phi^\top\frac{\partial^2 L_0(\vec{b}_0,\vec{\phi}_0)}{\partial \vec{\phi}\partial \vec{\phi}^\top}\vec{u}_\phi\right|(1+o(1))> c^*\tau_n
\]
for some constant $c^*>0$. 
This implies that
\[
|L_0(\vec{b},\vec{\phi})-L_0(\vec{b}_0,\vec{\phi}_0)| > c^*\tau_n\varepsilon^4.
\]
In following, we redefine $\varepsilon^4$ as $2^{-1}c^*\varepsilon^4$. 
Accordingly, we obtain 
\begin{eqnarray*}
&&P(\|\hat{\vec{b}}-\vec{b}_0\|^2+\|\vec{\phi}-\vec{\phi}_0\|^2>\varepsilon^2)\\
&&\leq P\left(\sup_{\eta(\vec{b},\vec{\phi})\leq \varepsilon^2}|L(\vec{b},\vec{\phi})-L_0(\vec{b},\vec{\phi})|\geq 2^{-1}\inf_{\eta(\vec{b},\vec{\phi})= \varepsilon^2}|L_0(\vec{b})-L_0(\vec{b}_0)| \right)\\
&&\leq P\left(\sup_{\eta(\vec{b},\vec{\phi})\leq \varepsilon^2}|L(\vec{b},\vec{\phi})-L_0(\vec{b},\vec{\phi})|\geq \varepsilon^4\tau_n \right).
\end{eqnarray*}

Then, our porpose is to show
\begin{eqnarray}
P\left(\sup_{\eta(\vec{b},\vec{\phi})\leq \varepsilon^2}|L(\vec{b},\vec{\phi})-L_0(\vec{b},\vec{\phi})|\geq \varepsilon^4\tau_n \right)\rightarrow 0. \label{ps2}
\end{eqnarray}
Again, we consider  $\vec{b}=\vec{b}_0+ \delta_n \vec{u}, \vec{u}\in\mathbb{R}^{K}$, where $\|\vec{u}\|\leq 1$. 
We then obtain
\begin{eqnarray*}
&&P\left(\sup_{\eta(\vec{b},\vec{\phi})\leq \varepsilon^2}|L(\vec{b},\vec{\phi})-L_0(\vec{b},\vec{\phi})|\geq \varepsilon^4\tau_n\right)\\
&&\leq P\left(|L(\vec{b}_0,\vec{\phi}_0)-L_0(\vec{b}_0,\vec{\phi}_0)|\geq 2^{-1} \varepsilon^4\tau_n\right)\\
&&+P\left(\sup_{\eta(\vec{b},\vec{\phi})\leq \varepsilon^2}|L(\vec{b},\vec{\phi})-L(\vec{b}_0,\vec{\phi}_0)-L_0(\vec{b},\vec{\phi})+L_0(\vec{b}_0,\vec{\phi}_0)|\geq 2^{-1}\varepsilon^4\tau_n\right)\\
&&\equiv J_1+J_2.
\end{eqnarray*}
We evaluate $J_1$. 
Define
\[
h(Y,\vec{X})=\left\{\exp[\vec{B}(\vec{X}^\top\vec{\theta}(\vec{\phi}_0))^\top\vec{b}_0]\log\left(\frac{Y}{w_n}\right)-\vec{B}(\vec{X}^\top\vec{\theta}(\vec{\phi}_0))^\top\vec{b}_0\right\}I(Y>w_n).
\]
From Lemma \ref{appspline}, we obtain 
\[
L(\vec{b}_0,\vec{\phi}_0)-L_0(\vec{b}_0,\vec{\phi}_0)=\frac{1}{n}\sum_{i=1}^n h(Y_i, \vec{X}_i)-E[h(Y_i, \vec{X}_i)] +\frac{\lambda}{2}\int_a^b \{\alpha_0^{(m)}(x)\}^2dx(1+o(1)).
\]
Under (C6), we have $\lambda/\tau_n =O(K^{-2m})=o(1)$.
Therefore, to show ${\cal J}_1\rightarrow 0$, it is sufficient to derive
\[
P\left(\left|n^{-1}\sum_{i=1}^n h(Y_i, \vec{X}_i)-E[h(Y_i, \vec{X}_i)] \right|> \varepsilon^4\tau_n \right) \rightarrow 0.
\]
Since $\exp[\vec{B}(\vec{X}^\top\vec{\theta}(\vec{\phi}_0))^\top\vec{b}_0]=\exp[\alpha(\vec{X}^\top\vec{\theta}_0))](1+o(1))$ and $\exp[\alpha(\vec{X}_i^\top\vec{\theta}_0)]\log(Y_i/w_n)$ is asymptotically distributed as standard exponential distribution under $Y_i>w_n$, $V[h(Y_i, \vec{X}_i)] \leq c^* \tau_n$ for some constant $c^*>0$. 
Therefore, Chebyshev's inequality and (C4) yield that 
\[
{\cal J}_1 = P\left(\left|n^{-1}\sum_{i=1}^n h(Y_i, \vec{X}_i)-E[h(Y_i,\vec{X}_i)] \right|> \varepsilon^4\tau_n \right)
\leq \frac{c^*}{n\tau_n\varepsilon^8}\rightarrow 0.
\]

Next, we focus on $J_2$. 
The Taylor expansion yields that
\begin{eqnarray*}
L(\vec{b},\vec{\phi})-L(\vec{b}_0,\vec{\phi}_0)
&=&
\varepsilon^2\left\{\frac{\partial L(\vec{b}_0,\vec{\phi}_0)}{\partial \vec{b}^\top}\vec{u}_b + \frac{\partial L(\vec{b}_0,\vec{\phi}_0)}{\partial \vec{\phi}^\top}\vec{u}_\phi + \lambda\vec{u}_b^\top \Delta_{m,K} \vec{b}_0\right\}(1+o(1))
\end{eqnarray*}
and 
\[
L_0(\vec{b},\vec{\phi})-L_0(\vec{b}_0,\vec{\phi}_0)
=
\varepsilon^2\left\{\frac{\partial L_0(\vec{b}_0,\vec{\phi}_0)}{\partial \vec{b}^\top}\vec{u}_b + \frac{\partial L_0(\vec{b}_0,\vec{\phi}_0)}{\partial \vec{\phi}^\top}\vec{u}_\phi\right\}(1+o(1)).
\]
By a similar argument as in the proof of Lemma \ref{gradient}, we obtain 
\[
\lambda\vec{u}_b^\top\Delta_{m,K} \vec{b}_0 =\lambda\vec{u}_b^\top D_{m,K}^\top \int \vec{B}^{[d-m]}(x) \alpha_0^{(m)}(x) dx =O(\lambda K^m).
\]
Under (C6), we have $O(\lambda K^m) = O(\tau_n K^{-m})=o(\tau_n)$, and hence the part $\lambda\vec{u}_b^\top \Delta_{m,K} \vec{b}_0$ is smaller than $\varepsilon^2\tau_n$. 
Thus, the remaining proof is to show
\begin{eqnarray*}
P\left(\sup_{\|\vec{u}\|^2< 1}\left|\left\{\frac{\partial L(\vec{b}_0,\vec{\phi}_0)}{\partial \vec{b}^\top}\vec{u}_b + \frac{\partial L(\vec{b}_0,\vec{\phi}_0)}{\partial \vec{\phi}^\top}\vec{u}_\phi \right\}  -\left\{\frac{\partial L_0(\vec{b}_0,\vec{\phi}_0)}{\partial \vec{b}^\top}\vec{u}_b + \frac{\partial L_0(\vec{b}_0,\vec{\phi}_0)}{\partial \vec{\phi}^\top}\vec{u}_\phi\right\}\right|\geq \tau_n\varepsilon^2\right)\rightarrow 0.
\end{eqnarray*}
Since 
\begin{eqnarray*}
&&P\left(\sup_{\|\vec{u}\|^2< 1}\left|\left\{\frac{\partial L(\vec{b}_0,\vec{\phi}_0)}{\partial \vec{b}^\top}\vec{u}_b + \frac{\partial L(\vec{b}_0,\vec{\phi}_0)}{\partial \vec{\phi}^\top}\vec{u}_\phi \right\}  -\left\{\frac{\partial L_0(\vec{b}_0,\vec{\phi}_0)}{\partial \vec{b}^\top}\vec{u}_b + \frac{\partial L_0(\vec{b}_0,\vec{\phi}_0)}{\partial \vec{\phi}^\top}\vec{u}_\phi\right\}\right|\geq \tau_n\varepsilon^2\right)\\
&&
\leq P\left(\sup_{\|\vec{u}_b\|^2< 1}\left|\frac{\partial (L(\vec{b}_0,\vec{\phi}_0)- L_0(\vec{b}_0,\vec{\phi}_0))}{\partial \vec{b}^\top}\vec{u}_b \right|\geq 2^{-1}\tau_n\varepsilon^2\right)\\
&&\quad +
P\left(\sup_{\|\vec{u}_\phi\|^2< 1}\left|\frac{\partial (L(\vec{b}_0,\vec{\phi}_0)- L_0(\vec{b}_0,\vec{\phi}_0))}{\partial \vec{\phi}^\top}\vec{u}_\phi \right|\geq 2^{-1}\tau_n\varepsilon^2\right)\\
&&\equiv {\cal J}_{21} + {\cal J}_{22}.
\end{eqnarray*}
From now on, we only show ${\cal J}_{21}\rightarrow 0$, but the proof of  ${\cal J}_{22}\rightarrow 0$ is similar. 

Let $E_i= \exp[\alpha_0(\vec{X}^\top\vec{\theta}_0)]\log(Y_i/w_n)$. 
Then, under $Y_i>w_n$, $E_i$ is approximately distributed as standard exponential distribution. 
From Lemma \ref{appspline} and proof of Lemma \ref{gradient}, we obtain 
\begin{eqnarray*}
\frac{\partial L(\vec{b}_0,\vec{\phi}_0)}{\partial \vec{b}^\top} \vec{u}_b &=& 
\frac{1}{n}\sum_{i=1}^n \left\{\exp[\vec{B}(\vec{X}_{i}^\top\vec{\theta}_0)^\top\vec{b}_0]\log\left(\frac{Y_i}{w_n}\right)-1 \right\}\vec{B}(\vec{X}_{i}^\top\vec{\theta}_0)^\top\vec{u}_bI(Y_i>w_n) \\
&=&
\frac{1}{n}\sum_{i=1}^n (E_i-1) \vec{B}(\vec{X}_{i}^\top\vec{\theta}_0)^\top\vec{u}_bI(Y_i>w_n)  + o_P(\tau_n).
\end{eqnarray*}
and 
\[
\frac{\partial L_0(\vec{b}_0,\vec{\phi}_0)}{\partial \vec{b}^\top} \vec{u}_b
= E[P(Y>w_n|\vec{X})r_n(X)B(\vec{X}^\top\vec{\theta}_0)^\top\vec{u}_b]=o(\tau_n). 
\]
Define the event ${\cal M}=\{\max_i E_i \leq \log(n)/\varepsilon^2\}$.  
We then have 
\[
{\cal J}_{21}
\leq P\left(\sup_{\|\vec{u}_b\|^2< 1}\left|\frac{1}{n}\sum_{i=1}^n (E_i-1) \vec{B}(\vec{X}_{i}^\top\vec{\theta}_0)^\top\vec{u}_bI(Y_i>w_n)\right| >\tau_n\varepsilon^2|{\cal M} \right)P({\cal M}) +P({\cal M}^c). 
\]
Since $P({\cal M})=(1-e^{-\log(n)/\varepsilon^2})^n$, we obtain $P({\cal M}^c) = 1-(1-e^{-\log(n)/\varepsilon^2})^n\rightarrow 0$.
Thus, the purpose is to show 
\[
P\left(\sup_{\|\vec{u}_b\|^2< 1}\left|\frac{1}{n}\sum_{i=1}^n (E_i-1) \vec{B}(\vec{X}_{i}^\top\vec{\theta}_0)^\top\vec{u}_bI(Y_i>w_n)\right| >\tau_n\varepsilon^2|{\cal M} \right)
\rightarrow 0.
\]
Let ${\cal U}=\{\vec{u}\in\mathbb{R}^{K} : \|\vec{u}\|< 1\}$ be the vector space and ${\cal U}_1,\ldots,{\cal U}_N$ be a covering of ${\cal U}$ with the diameter $R_n=C/(4n^\nu)$ for some constant $C>0$ and $\nu>0$. 
That is, ${\cal U}\subseteq \cup_{i=1}^N{\cal U}_i$. 
Then, Lemma 2.5 of van de Geer (2000) yields that it is sufficient to set $N\leq C (n^\nu )^{K}$. 
Let $\vec{u}_{j,b}\in{\cal U}_j , j=1,\ldots,N$. 
Then, for any $\vec{u}\in{\cal U}_j$, $\|\vec{u}-\vec{u}_{j,b}\|\leq R_n$. 
Therefore, we have 
\begin{eqnarray*}
&&\sup_{\|\vec{u}_b\|^2< 1}\left|\frac{1}{n}\sum_{i=1}^n (E_i-1) \vec{B}(\vec{X}_{i}^\top\vec{\theta}_0)^\top\vec{u}_bI(Y_i>w_n)\right| \\
&&
\leq\max_{1\leq j\leq N} \left|\frac{1}{n}\sum_{i=1}^n (E_i-1) \vec{B}(\vec{X}_{i}^\top\vec{\theta}_0)^\top\vec{u}_{j,b}I(Y_i>w_n)\right| \\
&&\quad + \max_{1\leq j\leq N}\sup_{\vec{u}_{b}\in{\cal U}_j} \left|\frac{1}{n}\sum_{i=1}^n (E_i-1) \vec{B}(\vec{X}_{i}^\top\vec{\theta}_0)^\top(\vec{u}_b-\vec{u}_{j,b})I(Y_i>w_n)\right|.
\end{eqnarray*}
Since the $B$-spline bases are non-negative and bounded functions and $n^{-1}\sum_{i=1}^n I(Y_i>w_n)=\tau_n(1+o(1))$, on the event ${\cal M}$, we obtain 
\begin{eqnarray*}
&&\sup_{\vec{u}_{b}\in{\cal U}_j}  \left|\frac{1}{n}\sum_{i=1}^n (E_i-1) \vec{B}(\vec{X}_{i}^\top\vec{\theta}_0)^\top(\vec{u}_b-\vec{u}_{j,b})I(Y_i>w_n)\right|\\
&&\leq \left(\sup_{z\in[a,b], \|\vec{v}\|=1} \{\vec{v}^\top \vec{B}(z)\}^2\right) |\log(n)/\varepsilon^2-1|  \tau_n\sup_{\vec{u}_{b}\in{\cal U}_j}  \|\vec{u}_b-\vec{u}_{j,b}\|
 \\
&&=  O_P(\tau_n\log(n)/n^\nu)\\
&&=o_P(\tau_n).
\end{eqnarray*}
Thus, we have 
\begin{eqnarray*}
&&\sup_{\|\vec{u}_b\|^2< 1}\left|\frac{1}{n}\sum_{i=1}^n (E_i-1) \vec{B}(\vec{X}_{i}^\top\vec{\theta}_0)^\top\vec{u}_bI(Y_i>w_n)\right| \\
&&
\leq\max_{1\leq j\leq N} \left|\frac{1}{n}\sum_{i=1}^n (E_i-1) \vec{B}(\vec{X}_{i}^\top\vec{\theta}_0)^\top\vec{u}_{j,b}I(Y_i>w_n)\right| +o_P(\tau_n).
\end{eqnarray*}
Lastly, we aim to derive 
\begin{eqnarray*}
P\left(\max_{1\leq j\leq N} \left|\frac{1}{n}\sum_{i=1}^n (E_i-1) \vec{B}(\vec{X}_{i}^\top\vec{\theta}_0)^\top\vec{u}_{j,b}I(Y_i>w_n)\right| >\tau_n\varepsilon^2|{\cal M} \right)
\rightarrow 0.
\end{eqnarray*}
We first obtain 
\begin{eqnarray*}
&&P\left(\max_{1\leq j\leq N} \left|\frac{1}{n}\sum_{i=1}^n (E_i-1) \vec{B}(\vec{X}_{i}^\top\vec{\theta}_0)^\top\vec{u}_{j,b}I(Y_i>w_n)\right| >\tau_n\varepsilon^2|{\cal M} \right)\\
&&\leq\sum_{j=1}^N P\left(\left|\frac{1}{n\tau_n}\sum_{i=1}^n (E_i-1) \vec{B}(\vec{X}_{i}^\top\vec{\theta}_0)^\top\vec{u}_{j,b}I(Y_i>w_n)\right| >\varepsilon^2|{\cal M} \right).
\end{eqnarray*}
On the event ${\cal M}$, it easy to find $(n\tau_n)^{-1}|(E_i-1) \vec{B}(\vec{X}_{i}^\top\vec{\theta}_0)^\top\vec{u}_{j,b}I(Y_i>w_n)|\leq C_1 \log n/(n\tau_n)$ for some constant $C_1>0$. 
Next, similar to proof of Lemma \ref{EigenHessian}, we have 
\[
V[(n\tau_n)^{-1}(E_i-1) \vec{B}(\vec{X}_{i}^\top\vec{\theta}_0)^\top\vec{u}_{j,b}I(Y_i>w_n)] \leq C_2 \frac{\log n}{n^2\tau_n}
\]
for some constant $C_2>0$. 
Therefore, Bernstein's inequality yields that 
\begin{eqnarray*}
P\left(\left|\frac{1}{n}\sum_{i=1}^n (E_i-1) \vec{B}(\vec{X}_{i}^\top\vec{\theta}_0)^\top\vec{u}_{j,b}I(Y_i>w_n)\right| >\tau_n\varepsilon^2|{\cal M} \right)
\leq C^* \exp\left[-C^*\varepsilon^4 n\tau_n/\log n \right]
\end{eqnarray*}
for some constant $C^*>0$. 
Therefore, under (C5), for some constants $C_0,C_1,C_2>0$, 
\begin{eqnarray*}
&&\sum_{j=1}^N P\left(\left|\frac{1}{n}\sum_{i=1}^n (E_i-1) \vec{B}(\vec{X}_{i}^\top\vec{\theta}_0)^\top\vec{u}_{j,b}I(Y_i>w_n)\right| >\tau_n\varepsilon^2|{\cal M} \right)\\
&&\leq C_0\exp[- C_1 \varepsilon^4 n\tau_n/\log n + C_2 K\log n ]\\
&&\rightarrow 0.
\end{eqnarray*}
Consequently, 
$
P(\|\hat{\vec{b}}-\vec{b}_0\|^2 + \|\hat{\vec{\phi}}-\vec{\phi}_0\|^2>\varepsilon^2)\rightarrow 0
$
was proven.
\end{proof}

\subsection*{Appendix C: Proof of Theorems}

\begin{proof}[Proof of Theorem \ref{RateParameter}]

From Lemma \ref{Consistency}, we have $\|\hat{\vec{b}}-\vec{b}_0\|\stackrel{P}{\rightarrow} 0$ and $\|\hat{\vec{\phi}}-\vec{\phi}_0\|\stackrel{P}{\rightarrow} 0$. 
Therefore, from te Taylors expansion of first derivative of penalized log-likehood function, we obtain 
\[
\left[
\begin{array}{c}
\hat{\vec{b}}-\vec{b}_0\\
\hat{\vec{\phi}}-\vec{\phi}_0
\end{array}
\right]
=\Sigma^{-1} 
\left[
\begin{array}{c}
\frac{\partial \ell_n(\vec{b}_0,\vec{\phi}_0)}
{\partial \vec{b}}\\
\frac{\partial \ell_n(\vec{b}_0,\vec{\phi}_0)}
{\partial \vec{\phi}}
\end{array}
\right](1+o(1)).
\]
From the property of inverse of block matrix, we obtain 
\[
\Sigma^{-1}=
\left[
\begin{array}{cc}
\Sigma_{b,b}& \Sigma_{b,\phi}\\
\Sigma_{\phi,b}&\Sigma_{\phi,\phi}
\end{array}
\right]^{-1}
=
\left[
\begin{array}{cc}
\Sigma_{b,b}^{-1}+\Sigma_{b,b}^{-1}\Sigma_{b,\phi} S_{\phi,\phi}\Sigma_{\phi,b}\Sigma_{b,b}^{-1}, & -\Sigma_{b,b}^{-1}\Sigma_{b,\phi}S_{\phi,\phi}\\
-S_{\phi,\phi}\Sigma_{\phi,b}\Sigma_{b,b}^{-1}, & S_{\phi,\phi}
\end{array}
\right]
\]
with $S_{\phi,\phi}= (\Sigma_{\phi,\phi}-\Sigma_{\phi,b}\Sigma_{b,b}^{-1}\Sigma_{b,\phi})^{-1}$. 
Similar to the proof of Lemma \ref{EigenHessian}, for any non-zero vector $\vec{v}\in\mathbb{R}^K$ with $\|\vec{v}\| <C, C>0$, all elements of $\Sigma_{\phi,b}\vec{v}$ and $\vec{v}^\top\Sigma_{b,\phi}$ has an order $O(\tau_n)$.  
Meanwhile, all elements of $\Sigma_{b,b}$ have $O(\tau_n)$ from Lemma \ref{EigenHessian}. 
In addition, since $\Sigma_{b,b}$ is band matrix, from the property of inverse of band matrix in Theorem 2.2 of Demko (1977), the order of each element of $\Sigma_{b,b}^{-1}\vec{v}$ is bounded by $O(\tau_n^{-1})$. 
Therefore, each element of $\Sigma_{b,b}^{-1}\Sigma_{b,\phi}$ has an order $O(1)$ and $\Sigma_{\phi,b}\Sigma_{b,b}^{-1}\Sigma_{b,\phi}=O(\tau_n)$. 
This yields that $S_{\phi,\phi}=O(\tau_n^{-1})$. 
Similarly, we obtain 
\[
\Sigma_{b,b}^{-1}+\Sigma_{b,b}^{-1}\Sigma_{b,\phi} S_{\phi,\phi}\Sigma_{\phi,b}\Sigma_{b,b}^{-1}=O(\tau_n^{-1})
\]
and $S_{\phi,\phi}\Sigma_{\phi,b}\Sigma_{b,b}^{-1} = O(\tau_n^{-1})$. 
We note that $\rho_{max}(\Sigma_{b,b}^{-2})=O(K\tau_n^{-2})$ even if $\rho_{max}(\Sigma_{b,b}^{-1})=O(\tau_n^{-1})$ since $\Sigma_{b,b}$ is $K$-square matrix. 
Thus, we have 
\[
\|\hat{\vec{b}}-\vec{b}_0\|^2
\leq  O(K\tau_n^{-2})\left\{ \left\| \frac{\partial \ell_n(\vec{b}_0,\vec{\phi}_0)}
{\partial \vec{b}}\right\|^2+ \left\|\frac{\partial \ell_n(\vec{b}_0,\vec{\phi}_0)}
{\partial \vec{\phi}}\right\|^2\right\}.
\]
Furthermore, from the property of Fisher information matrix and Lemmas \ref{gradient}--\ref{hessian}, we have 
\[
E\left[\left\| \frac{\partial \ell_n(\vec{b}_0,\vec{\phi}_0)}
{\partial \vec{b}}\right\|^2\right]\leq O\left(\frac{\tau_n}{n}\right) + O(\tau_n^{2+ 2\beta_{inf}} K^{-1}) + O(\tau_n\lambda K^{-1})
\]
and 
\[
E\left[\left\|\frac{\partial \ell_n(\vec{b}_0,\vec{\phi}_0)}
{\partial \vec{\phi}}\right\|^2\right] \leq O\left(\frac{\tau_n}{n}\right) + O(\tau_n^{2+ 2\beta_{inf}}).
\]
Since $K=O((\lambda/\tau_n)^{-1/(2m)})$ by (C6), we have 
\[
E[\|\hat{\vec{b}}-\vec{b}_0\|^2] \leq O\left(\frac{1}{n\tau_n}\left(\frac{\lambda}{\tau_n}\right)^{-1/(2m)}\right)+ O(\tau_n^{2\beta_{\inf}}) + O(\lambda/\tau_n).
\]
Similarly, we can obtain 
\[
E[\|\hat{\vec{\phi}}-\vec{\phi}_0\|^2] \leq O\left(\frac{1}{n\tau_n}\right)+ O(\tau_n^{2\beta_{\inf}}).
\]
\end{proof}

\begin{proof}[Proof of Theorem \ref{RateEstimator}]

We remember $\hat{\vec{\theta}}=\vec{\theta}(\hat{\vec{\phi}})$ and $\vec{\theta}_0=\vec{\theta}(\vec{\phi}_0)$.
From Lemma \ref{Consistency} and the Taylor expansion, we have
\begin{eqnarray}
\hat{\alpha}(\vec{X}^\top \hat{\vec{\theta}})
&=&
\vec{B}(\vec{X}^\top \hat{\vec{\theta}})^\top\hat{\vec{b}} \nonumber\\
&=&\vec{B}(\vec{X}^\top \vec{\theta}_0)^\top\vec{b}_0 +\vec{B}(\vec{X}^\top \vec{\theta}_0)^\top(\hat{\vec{b}}-\vec{b}_0)(1+o_P(1)) \nonumber\\
&&+ \alpha_1^{(1)}(\vec{X}^\top\vec{\theta}_0)\vec{X}^\top J_1(\vec{\phi}_0)(\hat{\vec{\phi}}-\vec{\phi}_0)(1+o_P(1)). \label{Cross}
\end{eqnarray}
This and Lemma \ref{appspline} yield that
\begin{eqnarray*}
\hat{\alpha}(\vec{X}^\top \hat{\vec{\theta}})-\alpha_0(\vec{X}^\top\vec{\theta}_0)
&=&
\vec{B}(\vec{X}^\top \vec{\theta}_0)^\top(\hat{\vec{b}}-\vec{b}_0)(1+o_P(1))\\
&&+ \alpha_1^{(1)}(\vec{X}^\top\vec{\theta}_0)\vec{X}^\top J_1(\vec{\phi}_0)(\hat{\vec{\phi}}-\vec{\phi}_0)(1+o_P(1)) + O(K^{-q}).
\end{eqnarray*}
From the proof of Lemma \ref{EigenHessian}, we have $\rho_{max}(E[\vec{B}(\vec{X}^\top\vec{\theta}_0)\vec{B}(\vec{X}^\top\vec{\theta}_0)^\top])\leq C$ for some constant $C>0$. 
This implies that 
\[
E\left[\left\{\vec{B}(\vec{X}^\top \hat{\vec{\theta}})^\top(\hat{\vec{b}}-\vec{b}_0)\right\}^2\right]
\leq C E[\|\hat{\vec{b}}-\vec{b}_0\|^2].
\]
Meanwhile, the domain of $\vec{X}$ is compact, we have 
\[
E\left[\left\{\alpha_1^{(1)}(\vec{X}^\top\vec{\theta}_0)\vec{X}^\top J_1(\vec{\phi}_0)(\hat{\vec{\phi}}-\vec{\phi}_0)\right\}^2\right]
\leq \tilde{C} E[\|\hat{\vec{\phi}}-\vec{\phi}_0\|^2]
\]
for some constant $\tilde{C}>0$. 
After applying Cauchy–Schwarz inequality to (\ref{Cross}), this theorem can be proven.

\end{proof}

\noindent{\bf ACKNOWLEDGEMENTS}

The authors are grateful to the Associate Editor and the anonymous referees for their valuable comments and suggestions, which have led to important improvements in the paper. 
This research was partially financially supported by the JSPS KAKENHI (Grant Nos. 22K11935 and 23K28043).
We would like to thank FASTEKJAPAN(www.fastekjapan.com) for English language editing.

\vspace{3mm}

\noindent{\bf DATA AVAILABILITY STATEMENT}
The data which support the findings of this study are available from the corresponding author upon reasonable request.

\def\bibname{References}

\end{document}